\documentclass[11pt,a4paper]{article}
\usepackage{mathrsfs}
\usepackage{amsmath,amssymb,amsthm}
\usepackage{graphicx,booktabs,float}
\usepackage[margin=1in]{geometry}
\usepackage[T1]{fontenc}
\usepackage{lmodern}
\usepackage[english]{babel}
\usepackage{enumitem}
\usepackage{mathtools}
\usepackage{hyperref}
\hypersetup{hidelinks,pdftitle={Nash Bounds for Entropy-Regularized Linear-Quadratic Games with Continuous Types},pdfauthor={Jisheng Liu}}

\newcommand{\bR}{\mathbb{R}}
\newcommand{\bE}{\mathbb{E}}
\newcommand{\bP}{\mathbb{P}}

\DeclareMathOperator{\tr}{tr}
\allowdisplaybreaks[3]
\newcommand{\dmN}{\bar{m}_N}
\newcommand{\barA}{\bar A}
\newcommand{\barC}{\bar C}

\newtheorem{theorem}{Theorem}[section]
\newtheorem{proposition}[theorem]{Proposition}
\newtheorem{corollary}[theorem]{Corollary}
\newtheorem{lemma}[theorem]{Lemma}
\newtheorem{assumption}[theorem]{Assumption}
\theoremstyle{remark}
\newtheorem{remark}[theorem]{Remark}

\title{Nash Bounds for Entropy-Regularized Linear-Quadratic Games with Continuous Types}
\author{Jisheng Liu\thanks{School of Mathematical Sciences, Fudan University, China. Correspondence: \href{mailto:jsliu23@m.fudan.edu.cn}{\texttt{jsliu23@m.fudan.edu.cn}}.}}
\date{}

\addto\captionsenglish{}
\makeatletter
\long\def\@makecaption#1#2{%
  \vskip\abovecaptionskip
  \noindent\textbf{#1}\enspace#2\par
  \vskip\belowcaptionskip
}
\makeatother

\begin{document}
\maketitle

\begin{abstract}
We derive finite-population Nash bounds for entropy-regularized
linear-quadratic games with continuously type-dependent coefficients.
The diffusion depends on each player's state and action and on the
empirical state mean. For a prescribed mean path, we construct
decentralized Gaussian feedback policies and bound the gain from a
unilateral deviation using all types and state histories. For any fixed
type profile, the bound separates a term inversely proportional to
population size from the squared discrepancy between the prescribed path
and the average of the players' mean responses to that path. At a mean
field fixed point and fixed entropy weight, independent type draws give
the inverse-population bound for all players simultaneously at any fixed
confidence level. The proof uses an exact finite-game value gap and
estimates of the state moments. For fixed types and mean input,
exploration adds a term linear in the entropy weight to the
strict-control Nash error. In a subclass with no mean-state costs and
nonvanishing action and mean loadings in the diffusion, the slope is of
inverse-population order.
\end{abstract}
\par

% introduction.tex
\section{Introduction}
\label{sec:intro}

In exploratory stochastic control, the controller chooses a distribution
over actions, and entropy regularization yields Gaussian feedback in
linear-quadratic problems \cite{WZZ20}. In our model, changing the entropy
weight changes the action variance while leaving the feedback mean fixed.
With control-dependent diffusion, this variance changes state fluctuations.
We ask how much a player in a heterogeneous finite population can gain
by departing from an exploratory policy designed for the limiting
population.

Mean field games approximate large stochastic games by
requiring individual optimal responses to reproduce the
population behavior against which they were computed
\cite{HMC06,LasryLions07}. In the linear-quadratic (LQ)
setting, Bensoussan et al.\ \cite{BensoussanSungYamYung16}
study the existence and uniqueness of mean field equilibria.
Huang and Zhou \cite{HuangZhou20} relate the fixed point
construction to the asymptotic solvability of finite-player
games.

Huang, Caines and Malham\'e \cite{HMC07} study
cost-coupled LQG games with nonuniform individual dynamics,
using the distribution of individual parameters to aggregate
their responses and construct decentralized approximate Nash
equilibria. Nguyen and Huang \cite{NguyenHuang12} allow
the minor players' dynamic coefficients to vary continuously
over a compact parameter set in a major--minor game. Their finite-game
bound retains the discrepancy between the empirical and limiting populations.
Hu, Huang and Nie \cite{HuHuangNie18} consider finitely many
types of minor players with heterogeneous input constraints and
diffusion depending on the action and the common state average.

Control-dependent diffusion introduces another feature:
actions affect both expected motion and fluctuations.
Xu and Zhang \cite{XuZhang20} and Si and Shi \cite{SiShi26}
study LQ mean field games with controlled diffusion.
For homogeneous LQ games with classical controls, Huang and
Yang \cite{HuangYang21} establish an $O(1/N)$-Nash bound with
controlled individual and common noises, within the asymptotic
solvability framework.

Entropy-regularized control brings exploration into this
setting. Policy optimization and the execution of exploratory
policies are studied in \cite{GRZ24} and \cite{SzpruchTZ24},
respectively. For mean field games, Guo, Xu and
Zariphopoulou \cite{GXZ22} analyze entropy regularization
and learning, with explicit policies in a homogeneous LQ model
with controlled diffusion.
Firoozi and Jaimungal \cite{FirooziJaimungal22} construct
Gaussian exploratory equilibria and establish a
finite-population Nash approximation for multivariate LQG
games with finitely many subpopulations and uncontrolled
diffusion.

We study scalar, finite-horizon games with independent noises and
coefficients that vary continuously with type. The player's own state,
the action and the empirical state mean enter the diffusion. For a
given population, we bound the largest gain from a deviation using
all states and types. The bound includes the squared discrepancy
between the prescribed mean and the average of the typewise response means.
At a mean field fixed point it gives Nash convergence for deterministic
populations whose empirical type laws converge to the limiting law.
For independently sampled types, at a fixed point and fixed temperature,
the Nash error is $O(1/N)$ both after averaging over opponent types
and simultaneously for all players at any fixed confidence level.
The averaged bound also applies to mean inputs away from a fixed point.
At fixed types and mean input, the gain splits into the strict-control
gain and a term linear in the exploration temperature.

\begin{samepage}
Section~\ref{sec:model-results} formulates the game and
information available to players. Section~\ref{sec:nce-solvability}
constructs the mean field policies. Section~\ref{sec:asy-Nash}
proves their finite-population gain bound.
Section~\ref{sec:temperature} studies temperature dependence,
and Section~\ref{sec:numerics} presents numerical illustrations.
Section~\ref{sec:conclusions} concludes.\par
\end{samepage}
% model.tex
\section{Game and Information}
\label{sec:model-results}
\label{sec:asy-formulation}

\paragraph{Notation.}
For a density $\pi$, write $\pi(da)=\pi(a)\,da$ and set
\[
\mathcal H(\pi)=-\int\pi\log\pi,\qquad
D_{\rm KL}(\pi\|\nu)=\int\pi\log(\pi/\nu).
\]
The symbol $\mathbf1$ denotes a vector of ones, and $e_i$ the
$i$th coordinate vector. We use Euclidean vector norms $|\cdot|$ and the induced matrix
norm $\|\cdot\|_{\rm op}$. For vectors and matrices,
$\|x\|_\infty=\max_j|x_j|$ and
$\|A\|_\infty=\max_j\sum_k|A_{jk}|$; for time functions it is
the supremum norm. Matrix inequalities refer to the positive
semidefinite order unless entrywise order is specified.
Constants $C$ may change between estimates; their dependence is specified in the relevant statements.

Fix a horizon $T>0$, a discount rate $\beta>0$ and an entropy
parameter (temperature) $\gamma>0$. The finite game is defined for
given types $\theta_1,\ldots,\theta_N$. When the types are sampled, they
are independent with law $F$, which need not have a density. The driving
processes $W^1,\ldots,W^N$ are independent standard Brownian motions,
independent of the sampled types.

\begin{assumption}
\label{ass:asy-LQ}
The type space $\Theta\subset\bR^d$ is compact. The maps
$A,\barA,B,C,\barC,D,\sigma,\allowbreak Q,\bar Q,R,H,\bar H$ are continuous on $\Theta$.
For every type, $B\ne0$, $\sigma>0$, $Q,R,H>0$ and
$\bar Q,\bar H\ge0$. All players start from the same state $x_0$.
\end{assumption}

We first specify the dynamics under strict controls. For a realized
type vector
$\theta^N=(\theta_1,\ldots,\theta_N)$, write $A_j=A(\theta_j)$,
and similarly for the other coefficients. Let $\alpha_t^j\in\bR$ be
player $j$'s action and write $\dmN=N^{-1}\sum_jX^j$.
The state equations are
\begin{equation}
\label{eq:ordinary-state}
dX_t^j=\bigl(A_jX_t^j+\barA_j\dmN(t)+B_j\alpha_t^j\bigr)dt
+\bigl(\sigma_j+C_jX_t^j+\barC_j\dmN+D_j\alpha_t^j\bigr)dW_t^j,\;
X_0^j=x_0.
\end{equation}
The bars identify population-mean loadings: $\barA$ acts in the
drift and $\barC$ in the diffusion. The coefficients $B$ and $D$
load the action in these two terms, and $\sigma$ is the diffusion
intercept. For fixed types, let $\bE_W$ denote expectation over the Brownian
motions.

Relaxation averages the drift and instantaneous covariance over actions
\cite[Section~3]{WZZ20}, \cite[Equation~(2.2)]{GXZ22}.
At a fixed state, the generator is affine in the action law, so
the attainable drift and covariance pairs form a convex set.
The independent Brownian motions give no mixed covariances.
For the density used at time $t$, define its mean and variance by
\[
u_t^j=\int a\,\pi_{j,t}(da),\qquad
v_t^j=\int(a-u_t^j)^2\pi_{j,t}(da).
\]
At a state $x$, write $\bar x_N=N^{-1}\sum_jx_j$. The averaged
instantaneous variance is
\begin{equation}
\label{eq:covariance-average}
\begin{aligned}
&\int(\sigma_j+C_jx_j+\barC_j\bar x_N+D_ja)^2\pi_{j,t}(da)\\
&\quad=(\sigma_j+C_jx_j+\barC_j\bar x_N)^2
+2(\sigma_j+C_jx_j+\barC_j\bar x_N)D_ju_t^j+D_j^2((u_t^j)^2+v_t^j)\\
&\quad=(\sigma_j+C_jx_j+\barC_j\bar x_N+D_ju_t^j)^2+D_j^2v_t^j.
\end{aligned}
\end{equation}
The relaxed dynamics are
\begin{equation}
\label{eq:finite-game-state}
\begin{aligned}
dX_t^j={}&\bigl(A_jX_t^j+\barA_j\dmN
+B_ju_t^j\bigr)dt+\sqrt{a_{j,t}^N}\,dW_t^j,\\
a_{j,t}^N={}&(\sigma_j+C_jX_t^j+\barC_j\dmN+D_ju_t^j)^2
+D_j^2v_t^j,\qquad X_0^j=x_0.
\end{aligned}
\end{equation}
When $D_j\ne0$, action densities with the same mean and different
variances can therefore give different state laws. A Dirac action gives the generator of
\eqref{eq:ordinary-state}; the entropy-regularized game below uses
densities with respect to Lebesgue measure.

Fix a candidate profile $\widehat\pi^N$ of own-state Markov policies,
to be constructed in Section~\ref{sec:asy-MFG}. When player $i$
uses $\pi_i$, the opponents retain $\widehat\pi_{-i}^N$, and all
states evolve under $(\pi_i,\widehat\pi_{-i}^N)$. Player $i$ minimizes
\begin{align}
\label{eq:finite-game-cost}
J_i^N(\pi_i,\widehat\pi_{-i}^N;\theta^N)
={}&\bE_W\Bigg[\int_0^T e^{-\beta t}\bigg(
\frac{Q_i}2(X_t^i)^2
+\frac{\bar Q_i}2\dmN(t)^2\notag\\
&\qquad\qquad+\frac{R_i}2\int_{\bR}a^2\pi_{i,t}(da)
-\gamma \mathcal H(\pi_{i,t})\bigg)dt\notag\\
&\qquad+e^{-\beta T}\bigg(
\frac{H_i}2(X_T^i)^2
+\frac{\bar H_i}2\dmN(T)^2\bigg)\Bigg].
\end{align}
Entropy rewards dispersion of the action density.
The Nash bounds below use this cost unless we specify strict controls.

Each decentralized feedback uses the player's current state, type
and a precomputed mean input. As in \cite[Theorem~6.1]{HuangYang21},
we test these policies against unilateral deviations with access
to all players' states. Here the deviator also knows the full type
vector. For the distinction between local and global observations,
see \cite[Section~II-A]{LZ08}.

The types remain fixed during the game. For each own type
$\zeta\in\Theta$, let $\bP_i^\zeta$ fix $\theta_i=\zeta$, retain
the independent opponent types with law $F$, and retain the
independent Brownian motions. For an integrable functional $Y$,
\begin{equation}
\label{eq:own-type-expectation}
\bE_i^\zeta[Y]
=\int_{\Theta^{N-1}}\bE_W\!\left[Y(\zeta,\vartheta^{-i},W)\right]
F^{\otimes(N-1)}(d\vartheta^{-i}).
\end{equation}
This parameter definition includes types of zero $F$-mass.
Since $J_i^N$ already includes $\bE_W$, its outer expectation
$\bE_i^\zeta[J_i^N]$ averages only the opponent types.

For fixed types, call a rule admissible if its closed-loop system
has a nonexplosive strong solution, is stochastically pathwise
unique, and satisfies
\[
\bE_W\int_0^T\int_{\bR}a^2\pi_{i,t}(da)dt<\infty,
\qquad
\bE_W\int_0^T|\mathcal H(\pi_{i,t})|dt<\infty.
\]
Pathwise uniqueness compares nonanticipative solutions with the
same initial data and Brownian inputs on a common filtered space.
Let $\mathcal D(\bR)$ be the probability densities with the Borel
structure inherited from $L^1(\bR)$, and let
$\mathcal G_t^N=\sigma(X_s^j:0\le s\le t,\ 1\le j\le N)$ be completed
under the law of the deviating system. For a revealed type vector,
\[
\mathcal A_i(\theta^N)
=\left\{\pi_i:\begin{array}{l}
\pi_i\text{ is }\mathcal D(\bR)\text{-valued and }
(\mathcal G_t^N)\text{-progressively measurable},\\
\pi_i\text{ is admissible for }\theta^N
\end{array}\right\}.
\]

At a revealed type vector, the full-information gain is the
supremum of the cost reduction from a unilateral deviation,
\begin{equation}
\label{eq:fixed-gain}
\varepsilon_{i,N}(\theta^N)
= J_i^N(\widehat\pi_i^N,\widehat\pi_{-i}^N;\theta^N)
-\inf_{\pi_i\in\mathcal A_i(\theta^N)}J_i^N(\pi_i,\widehat\pi_{-i}^N;\theta^N).
\end{equation}
For an admissible candidate this gain is nonnegative. Its average
over opponent types is
\[
\mathcal E_{i,N}(\zeta)=\bE_i^\zeta[\varepsilon_{i,N}(\theta^N)].
\]

A decentralized deviation based on the player's own type and
state history has the form
$\pi_{i,t}=\phi_i(t,X^i_{\cdot\wedge t})$, where
$X^i_{\cdot\wedge t}=(X^i_{s\wedge t})_{0\le s\le T}$ and
$\phi_i:[0,T]\times C([0,T];\bR)\to\mathcal D(\bR)$ is Borel.
For fixed own type $\zeta$, the same map is used for every opponent
type vector. Suppose this rule is admissible for every such vector,
with $\bE_i^\zeta\int_0^T[\int a^2\pi_{i,t}(da)+|\mathcal H(\pi_{i,t})|]dt<\infty$.
The rule therefore belongs to $\mathcal A_i(\theta^N)$ for each such vector, so
\begin{equation}
\label{eq:own-observation-bound}
\bE_i^\zeta\!\left[
J_i^N(\widehat\pi_i^N,\widehat\pi_{-i}^N;\theta^N)
-J_i^N(\pi_i,\widehat\pi_{-i}^N;\theta^N)\right]
\le\mathcal E_{i,N}(\zeta).
\end{equation}
% mean_field.tex
\section{Mean Field Policy}
\label{sec:nce-solvability}

For each prescribed mean path, we construct the policies whose
finite-game deviation gain is bounded in Section~\ref{sec:asy-Nash}.
The consistency condition identifies when the policies generate the prescribed mean.

\subsection{Control for a Prescribed Mean}
\label{sec:asy-MFG}

The Nash certainty equivalence (NCE) construction requires optimality
for a prescribed mean
and consistency of the mean generated by that policy
\cite[Section~3.1]{HMC06}. Fix $m\in C([0,T])$ and a type $\theta$.
Write $A_\theta=A(\theta)$, and similarly for other coefficients.
The prescribed-mean diffusion intercept is
$\sigma_\theta[m](t)=\sigma_\theta+\barC_\theta m(t)$.
Replacing $\dmN$ by $m$ in \eqref{eq:finite-game-state} gives
\begin{equation}
\label{eq:limiting-state}
\begin{aligned}
dX_t={}&(A_\theta X_t+\barA_\theta m+B_\theta u_t)dt
+\sqrt{(\sigma_\theta[m]+C_\theta X_t+D_\theta u_t)^2+D_\theta^2v_t}\,dW_t,
\end{aligned}
\end{equation}
where $u_t$ and $v_t$ are the action density's mean and variance.
The initial state is $x_0$. The type-$\theta$ criterion is
\begin{equation}
\label{eq:limiting-cost}
\begin{aligned}
J^\theta(\pi;m)=\bE\bigg[&\int_0^T e^{-\beta t}
\bigg(\frac{Q_\theta}2X_t^2+\frac{\bar Q_\theta}2m^2 +\frac{R_\theta}2\int a^2\pi_t(da)
-\gamma\mathcal H(\pi_t)\bigg)dt \\
&+e^{-\beta T}\bigg(\frac{H_\theta}2X_T^2
+\frac{\bar H_\theta}2m(T)^2\bigg)\bigg].
\end{aligned}
\end{equation}
The limiting admissible policies are Markov in $(t,X_t)$ and
satisfy the moment, entropy and solution conditions of Section~2.

Let $V^\theta(t,x;m)$ denote the value with discount measured from $t$.
We construct a candidate through the HJB equation
\begin{equation}
\label{eq:limiting-HJB}
\begin{aligned}
&-V_t^\theta+\beta V^\theta
={}(A_\theta x+\barA_\theta m)V_x^\theta
+\frac{Q_\theta}2x^2+\frac{\bar Q_\theta}2m^2 \\
&+\inf_\pi\Bigl\{\int\Bigl[
\frac{R_\theta}2a^2+B_\theta V_x^\theta a +\frac{V_{xx}^\theta}2(\sigma_\theta[m]+C_\theta x+D_\theta a)^2
\Bigr]\pi(da)-\gamma\mathcal H(\pi)\Bigr\}, \\
&V^\theta(T,x;m)={}\frac{H_\theta}2x^2
+\frac{\bar H_\theta}2m(T)^2.
\end{aligned}
\end{equation}
The quadratic ansatz and effective action weight are
\begin{equation}
\label{eq:quadratic-value}
V^\theta=\tfrac12P_\theta x^2+\eta_\theta x+\kappa_\theta,
\qquad \widehat R_\theta=R_\theta+D_\theta^2P_\theta.
\end{equation}
The HJB coefficients of $a^2$ and $a$ are $\widehat R_\theta/2$
and $B_\theta V_x^\theta+D_\theta P_\theta(\sigma_\theta[m]+C_\theta x)$.
For $\widehat R_\theta>0$, completing the square gives
\begin{equation}
\label{eq:mean-policy}
\begin{aligned}
u_\theta^m(t,x)&
=-\frac{(B_\theta+D_\theta C_\theta)P_\theta x+B_\theta\eta_\theta+D_\theta P_\theta\sigma_\theta[m]}
{\widehat R_\theta},\\
\pi_\theta^m(da\mid t,x)&=\mathcal N\!\left(u_\theta^m(t,x),\frac{\gamma}{\widehat R_\theta}\right)(da).
\end{aligned}
\end{equation}
Expanding the Gaussian log density gives
\begin{equation}
\label{eq:gibbs-identity}
\int\frac{\widehat R_\theta}{2}(a-u_\theta^m(t,x))^2\,\pi(da)-\gamma\mathcal H(\pi)
=-\frac{\gamma}{2}\log\frac{2\pi\gamma}{\widehat R_\theta}
+\gamma D_{\rm KL}(\pi\|\pi_\theta^m(\cdot\mid t,x)).
\end{equation}
The KL term is nonnegative and vanishes at $\pi_\theta^m$. The
infimum in \eqref{eq:limiting-HJB} is therefore
$P_\theta(\sigma_\theta[m]+C_\theta x)^2/2-\widehat R_\theta(u_\theta^m(t,x))^2/2
-\frac{\gamma}{2}\log(2\pi\gamma/\widehat R_\theta)$
\cite{WZZ20,GXZ22,GRZ24}.
The coefficients of $x^2$ and $x$ in the HJB satisfy
\[
\begin{aligned}
\tfrac12(-\dot P_\theta+\beta P_\theta)
&=(A_\theta+C_\theta^2/2)P_\theta+Q_\theta/2
-\frac{(B_\theta+D_\theta C_\theta)^2P_\theta^2}{2\widehat R_\theta},\\
-\dot\eta_\theta+\beta\eta_\theta
&=A_\theta\eta_\theta+\barA_\theta P_\theta m+C_\theta P_\theta\sigma_\theta[m]\\
&\quad-\frac{(B_\theta+D_\theta C_\theta)P_\theta(B_\theta\eta_\theta+D_\theta P_\theta\sigma_\theta[m])}{\widehat R_\theta}.
\end{aligned}
\]
The first equation gives
\begin{equation}
\label{eq:asy-Riccati}
\dot P_\theta=(\beta-2A_\theta-C_\theta^2)P_\theta
+\frac{(B_\theta+D_\theta C_\theta)^2P_\theta^2}{\widehat R_\theta}-Q_\theta,
\qquad P_\theta(T)=H_\theta.
\end{equation}
Write $u_\theta^m(t,x)=-K_\theta(t)x-k_\theta[m](t)$ with
\begin{equation}
\label{eq:policy-coefficients}
K_\theta=\frac{(B_\theta+D_\theta C_\theta)P_\theta}{\widehat R_\theta},\qquad
k_\theta[m]=\frac{B_\theta\eta_\theta
+D_\theta P_\theta\sigma_\theta[m]}{\widehat R_\theta}.
\end{equation}
With $A_\theta^{\rm cl}=A_\theta-B_\theta K_\theta$, the
linear and constant coefficients give
\begin{equation}
\label{eq:asy-eta-ODE}
\dot\eta_\theta=(\beta-A_\theta^{\rm cl})\eta_\theta
-P_\theta\bigl[\barA_\theta m+(C_\theta-D_\theta K_\theta)\sigma_\theta[m]\bigr],
\qquad \eta_\theta(T)=0.
\end{equation}
The terminal constant is $\kappa_\theta(T)=\bar H_\theta m(T)^2/2$, and
\begin{equation}
\label{eq:kappa}
\dot\kappa_\theta=\beta\kappa_\theta-\barA_\theta\eta_\theta m
-\tfrac12\bar Q_\theta m^2-\tfrac12P_\theta\sigma_\theta[m]^2
+\tfrac12\widehat R_\theta k_\theta[m]^2
+\tfrac\gamma2\log\frac{2\pi\gamma}{\widehat R_\theta}.
\end{equation}
Solve first for $P_\theta$, independently of $m$ and $\gamma$,
then for $\eta_\theta$ and $\kappa_\theta$ in that order.

\subsection{Verification}
\label{sec:limiting-verification}

In reversed time $p_\theta(\tau)=P_\theta(T-\tau)$ satisfies
\[
p_\theta'=Q_\theta+(2A_\theta+C_\theta^2-\beta)p_\theta
-\frac{(B_\theta+D_\theta C_\theta)^2p_\theta^2}{R_\theta+D_\theta^2p_\theta},
\qquad p_\theta(0)=H_\theta>0.
\]
At a first zero its left derivative would be nonpositive, whereas
the equation gives $p_\theta'=Q_\theta>0$. Hence $p_\theta>0$.
On this region $p_\theta'\le C(1+p_\theta)$, so
\[
0<p_\theta(\tau)\le(\sup_\theta H_\theta+C\tau)e^{C\tau}.
\]
This excludes explosion and gives $\widehat R_\theta\ge\inf_\theta R_\theta>0$.
On this bounded interval the Riccati field is uniformly
Lipschitz in $P$. Subtracting the terminal equations gives
\[
\begin{aligned}
|P_\theta(t)-P_\vartheta(t)|
&\le |H_\theta-H_\vartheta|
+C\int_t^T|P_\theta(s)-P_\vartheta(s)|ds\\
&\quad+CT\sum_{c\in\{A,B,C,D,Q,R\}}|c_\theta-c_\vartheta|.
\end{aligned}
\]
Gronwall and continuity of the coefficient maps imply continuity
of $P$ in type, uniformly in time.
The coefficients $P,\widehat R,K,A^{\rm cl}$ are therefore continuous in type.
The linear terminal equations determine continuous $\eta$ and
$\kappa$ for each prescribed $m$.

For the It\^o verification (cf. \cite[Lemma~3.2]{GRZ24}),
we first bound the state moments in \eqref{eq:limiting-state}
under an admissible policy $\pi$.
Its action mean and variance satisfy
\[
u_t^2\le\int a^2\pi_t(da),\qquad
(\sigma_\theta[m]+C_\theta X_t+D_\theta u_t)^2+D_\theta^2v_t
\le C_m\left(1+|X_t|^2+\int a^2\pi_t(da)\right).
\]
We stop the state and accumulated action moment at
\begin{equation}
\label{eq:scalar-stopping}
\tau_n=T\wedge\inf\left\{t:|X_t|\ge n\text{ or }
\int_0^t\int a^2\pi_s(da)ds\ge n\right\}.
\end{equation}
Nonexplosion and the finite action moment give $\tau_n\uparrow T$.
The integral state equation, Cauchy--Schwarz and Doob's inequality give
\[
\begin{aligned}
&\bE\sup_{s\le t}|X_{s\wedge\tau_n}|^2\le 3|x_0|^2
+3T\bE\int_0^{t\wedge\tau_n}|A_\theta X_s+\barA_\theta m+B_\theta u_s|^2ds\\
&\quad+12\bE\int_0^{t\wedge\tau_n}
\left[(\sigma_\theta[m]+C_\theta X_s+D_\theta u_s)^2+D_\theta^2v_s\right]ds\\
&\le C_m\left(1+|x_0|^2+\bE\int_0^T\int a^2\pi_s(da)ds\right)
+C_m\int_0^t\bE\sup_{r\le s}|X_{r\wedge\tau_n}|^2ds.
\end{aligned}
\]
Gronwall and Fatou yield
\begin{equation}
\label{eq:verification-moment}
\bE\sup_{t\le T}|X_t|^2
\le C_m\left(1+|x_0|^2+\bE\int_0^T\int a^2\pi_t(da)dt\right)<\infty.
\end{equation}

For the quadratic candidate,
$V_x^\theta=P_\theta X_t+\eta_\theta$ and $V_{xx}^\theta=P_\theta$.
Substituting these derivatives and the coefficient equations
into the drift plus running cost gives
\[
\begin{aligned}
&V_t^\theta-\beta V^\theta
+(A_\theta X_t+\barA_\theta m+B_\theta u_t)(P_\theta X_t+\eta_\theta)
+\frac{P_\theta}2\left[(\sigma_\theta[m]+C_\theta X_t+D_\theta u_t)^2+D_\theta^2v_t\right]\\
&\quad+\frac{Q_\theta}2X_t^2+\frac{\bar Q_\theta}2m^2
+\frac{R_\theta}2(u_t^2+v_t)-\gamma\mathcal H(\pi_t)\\
&=\frac{\widehat R_\theta}2\left[(u_t+K_\theta X_t+k_\theta[m])^2+v_t\right]
-\gamma\mathcal H(\pi_t)+\frac\gamma2\log\frac{2\pi\gamma}{\widehat R_\theta}\\
&=\gamma D_{\rm KL}(\pi_t\|\pi_\theta^m(\cdot\mid t,X_t)).
\end{aligned}
\]
The last equality is \eqref{eq:gibbs-identity}.
Let the stochastic integral in the discounted It\^o formula be
\[
\mathcal M_s=\int_0^s e^{-\beta t}(P_\theta X_t+\eta_\theta)
\sqrt{(\sigma_\theta[m]+C_\theta X_t+D_\theta u_t)^2+D_\theta^2v_t}\,dW_t.
\]
The Burkholder--Davis--Gundy (BDG) inequality and Cauchy--Schwarz give
\begin{equation}
\label{eq:verification-martingale-bound}
\begin{aligned}
\bE\sup_{s\le T}|\mathcal M_s|
&\le C_m\bE\left[
\left(1+\sup_{t\le T}|X_t|\right)
\left(\int_0^T[(\sigma_\theta[m]+C_\theta X_t+D_\theta u_t)^2+D_\theta^2v_t]dt\right)^{1/2}\right]\\
&\le C_m\left[\bE\left(1+\sup_{t\le T}|X_t|^2\right)\right]^{1/2}
\left[1+\bE\int_0^T\left(|X_t|^2+\int a^2\pi_t(da)\right)dt\right]^{1/2}<\infty.
\end{aligned}
\end{equation}
Here the first factor is finite by \eqref{eq:verification-moment},
and the second by the same estimate and the action-moment bound.
Thus $\mathcal M$ is a uniformly integrable martingale and
$\bE\mathcal M_T=0$. The endpoint and running cost are integrable
by quadratic growth, \eqref{eq:verification-moment} and the
absolute-entropy condition.
Integrating It\^o's formula gives the identity
\[
\begin{aligned}
&e^{-\beta T}V^\theta(T,X_T;m)-V^\theta(0,x_0;m)\\
&+\int_0^T e^{-\beta t}\left[
\frac{Q_\theta}2X_t^2+\frac{\bar Q_\theta}2m^2
+\frac{R_\theta}2\int a^2\pi_t(da)-\gamma\mathcal H(\pi_t)\right]dt\\
&=\gamma\int_0^T e^{-\beta t}
D_{\rm KL}(\pi_t\|\pi_\theta^m(\cdot\mid t,X_t))dt+\mathcal M_T.
\end{aligned}
\]
Taking expectations and using the terminal condition in
\eqref{eq:limiting-HJB} yields
\begin{equation}
\label{eq:limiting-verification-ineq}
J^\theta(\pi;m)-V^\theta(0,x_0;m)
=\gamma\bE\int_0^T e^{-\beta t}
D_{\rm KL}(\pi_t\|\pi_\theta^m(\cdot\mid t,X_t))dt\ge0.
\end{equation}
Under \eqref{eq:mean-policy} the drift is affine and the diffusion
is the norm of the affine two-vector
$(\sigma_\theta[m]+(C_\theta-D_\theta K_\theta)x-D_\theta k_\theta,
D_\theta\sqrt{\gamma/\widehat R_\theta})$.
It is globally Lipschitz in $x$ and has linear growth.
The unique strong solution has finite maximal second moment.
Moreover,
\[
\begin{aligned}
\int a^2\pi_\theta^m(da\mid t,X_t)
&=(K_\theta X_t+k_\theta[m])^2+\gamma/\widehat R_\theta
\le C_m(1+|X_t|^2)+C\gamma,\\
\mathcal H(\pi_\theta^m(\cdot\mid t,X_t))
&=\tfrac12\log(2\pi e\gamma/\widehat R_\theta).
\end{aligned}
\]
The entropy is bounded in time, so this policy is admissible.
Its KL term vanishes, proving optimality.

\subsection{Consistency}
\label{subsec:consistency}

For a fixed type, put $\mu_\theta[m](t)=\bE X_t$ under \eqref{eq:mean-policy}.
In \eqref{eq:limiting-state}, the stochastic integral has zero mean
by square integrability. The feedback mean is affine, so
$\bE[u_\theta^m(t,X_t)]=u_\theta^m(t,\mu_\theta[m](t))$.
Taking expectations in the drift gives
\begin{equation}
\label{eq:mu-theta}
\dot \mu_\theta=A_\theta^{\rm cl} \mu_\theta+\barA_\theta m-B_\theta k_\theta[m],
\qquad \mu_\theta(0)=x_0.
\end{equation}
For later estimates, variation of constants in
\eqref{eq:asy-eta-ODE} and \eqref{eq:mu-theta} gives
\begin{equation}
\label{eq:affine-coefficient-formulas}
\begin{aligned}
\eta_\theta(t)&=\int_t^T e^{-\int_t^s(\beta-A_\theta^{\rm cl}(r))dr}P_\theta(s)\bigl[\barA_\theta m(s)
+(C_\theta-D_\theta K_\theta(s))\sigma_\theta[m](s)\bigr]ds, \\
\mu_\theta[m](t)&=e^{\int_0^tA_\theta^{\rm cl}(r)dr}x_0
+\int_0^t e^{\int_s^tA_\theta^{\rm cl}(r)dr}
\bigl[\barA_\theta m(s)-B_\theta k_\theta[m](s)\bigr]ds.
\end{aligned}
\end{equation}
Thus for $\|m\|_\infty\le M$,
\begin{equation}
\label{eq:uniform-affine-coefficients}
\sup_\theta\bigl(\|\eta_\theta[m]\|_\infty+
\|k_\theta[m]\|_\infty+\|\mu_\theta[m]\|_\infty\bigr)\le C_M,
\end{equation}
with $C_M$ independent of $\gamma$. These functions are jointly
continuous in time and type. Define
\begin{equation}
\label{eq:NCE-map}
\mathscr T[m](t)=\int_\Theta \mu_\theta[m](t)F(d\theta).
\end{equation}
Joint continuity and \eqref{eq:uniform-affine-coefficients} make
$\mathscr T:C([0,T])\to C([0,T])$ well defined.
For any input $m$, the finite-game candidate is
\begin{equation}
\label{eq:asy-plugin}
\widehat\pi_i^N[m](da\mid t,x,\theta_i)=\pi_{\theta_i}^{m}(da\mid t,x).
\end{equation}
Each player uses its own state and type and the precomputed path $m$.
Neither the policy mean nor $\mathscr T$ depends on $\gamma$. The
Gaussian action variance $\gamma/\widehat R_\theta$ enters the
diffusion in \eqref{eq:finite-game-state}.

The discrepancy $\mathscr T[m]-m$ is the consistency residual.
A mean field fixed point is a path $m^*\in C([0,T])$ satisfying
$m^*=\mathscr T[m^*]$. At such a path the policies are optimal for
the mean they generate. Their state laws need not be Gaussian.

For two inputs $m_1,m_2$, denote the corresponding differences by
$\delta m,\delta\eta_\theta,\delta\mu_\theta$. In the backward
and forward equations, the loadings of $\delta m$ are different.
Write
\[
b_\theta^\eta=\barA_\theta+(C_\theta-D_\theta K_\theta)\barC_\theta,
\qquad
b_\theta^\mu=\barA_\theta-\frac{B_\theta D_\theta P_\theta\barC_\theta}{\widehat R_\theta}.
\]
Subtracting \eqref{eq:asy-eta-ODE} and \eqref{eq:mu-theta}, respectively, gives
\[
\begin{aligned}
\delta\dot\eta_\theta&=(\beta-A_\theta^{\rm cl})\delta\eta_\theta
-P_\theta b_\theta^\eta\delta m, & \delta\eta_\theta(T)&=0,\\
\delta\dot\mu_\theta&=A_\theta^{\rm cl}\delta\mu_\theta
+b_\theta^\mu\delta m-\frac{B_\theta^2}{\widehat R_\theta}\delta\eta_\theta,
&\delta\mu_\theta(0)&=0.
\end{aligned}
\]
With suprema over time and type for coefficient families, variation
of constants gives
\[
\sup_\theta\|\delta\eta_\theta\|_\infty
\le T e^{T\|\beta-A^{\rm cl}\|_\infty}\|P\|_\infty
\|b^\eta\|_\infty\|\delta m\|_\infty.
\]
Using this bound in the forward equation yields
$\|\mathscr T[m_1]-\mathscr T[m_2]\|_\infty\le L_T\|\delta m\|_\infty$, where
\begin{equation}
\label{eq:contraction-bound}
L_T=e^{T\|A^{\rm cl}\|_\infty}
\left[T\|b^\mu\|_\infty+
T^2\left\|\frac{B^2}{\widehat R}\right\|_\infty
\|P\|_\infty\|b^\eta\|_\infty e^{T\|\beta-A^{\rm cl}\|_\infty}\right].
\end{equation}
Thus $L_T<1$ for small $T$, or for sufficiently small mean loadings
$\barA,\barC$ at fixed remaining parameters.

\begin{proposition}
\label{prop:consistency}
If $\mathscr T$ is a contraction, it has a unique fixed point.
For a one-point type law, a sufficient condition for a unique fixed point on
$[0,T]$ is
\begin{equation}
\label{eq:homogeneous-sign}
Q+\frac{CP R}{\widehat R}(C+\barC)\ge0\qquad\text{on }[0,T].
\end{equation}
This condition holds for every finite $T$ if $C(C+\barC)\ge0$.
\end{proposition}
\begin{proof}
Banach's theorem gives the contraction assertion. For a one-point
type law, suppress the type and impose $\mu[m]=m$. Set $y=-Pm-\eta$.
Equations \eqref{eq:asy-Riccati}, \eqref{eq:asy-eta-ODE}
and \eqref{eq:mu-theta} give
\begin{equation}
\label{eq:homogeneous-system}
\begin{aligned}
\dot m&=\left[A+\barA-\frac{BDP(C+\barC)}{\widehat R}\right]m
+\frac{B^2}{\widehat R}y-\frac{BDP\sigma}{\widehat R},\\
\dot y&=\left[Q+\frac{CP R}{\widehat R}(C+\barC)\right]m
+\left[\beta-A+\frac{BDCP}{\widehat R}\right]y
+\frac{CP R\sigma}{\widehat R},\\
m(0)&=x_0,\qquad y(T)=-Hm(T).
\end{aligned}
\end{equation}
Conversely, from a solution of this system, $\eta=-y-Pm$ solves
\eqref{eq:asy-eta-ODE} and $m=\mu[m]$ by \eqref{eq:mu-theta}.

Under \eqref{eq:homogeneous-sign}, the coefficient matrix of
\eqref{eq:homogeneous-system} has nonnegative off-diagonal entries.
Its transition matrix $\Phi(t,s)$ is entrywise nonnegative, and
$\Phi_{22}(T,0)>0$ by the integrating factor for its second row.
The terminal residual $y(T)+Hm(T)$ is an affine function of $y(0)$
with slope $\Phi_{22}(T,0)+H\Phi_{12}(T,0)>0$
(cf.~\cite[Proposition~10]{HuangZhou20}). Thus there is exactly one
initial value $y(0)$ for which the residual vanishes.
\end{proof}
% finite_response.tex
\section{Finite-Population Performance}
\label{sec:asy-Nash}

For the candidate \eqref{eq:asy-plugin} constructed from $m$, write
$\varepsilon_{i,N}(\theta^N;m)$ for the gain in \eqref{eq:fixed-gain}
and $\mathcal E_{i,N}(\zeta;m)$ for its average over opponent types.

\begin{theorem}
\label{thm:main}
Under Assumption~\ref{ass:asy-LQ}, fix $M<\infty$. There is $C_M$,
depending only on $M,T,\beta,x_0$ and the coefficient bounds, such
that, for all $N\ge2$, $i\in\{1,\ldots,N\}$, $\zeta\in\Theta$,
$\gamma>0$ and continuous
inputs with $\|m\|_\infty\le M$,
\begin{equation}
\label{eq:main-residual}
0\le\mathcal E_{i,N}(\zeta;m)
\le C_M\left(\frac{1+\gamma}{N}
+\|\mathscr T[m]-m\|_{L^2(0,T)}^2\right).
\end{equation}
At any mean field fixed point $m^*$, the decentralized profile
has a $C(1+\gamma)/N$ Nash bound against full-information
deviations, averaged over opponent types and uniform in the
player's own type. By \eqref{eq:own-observation-bound}, the same
bound holds for decentralized deviations based on the player's
own type and state history.
\end{theorem}

The first corollary concerns a population with given types. The
second concerns a population sampled from $F$.

\begin{corollary}
\label{cor:fixed-types}
Under Assumption~\ref{ass:asy-LQ}, let $m$ be continuous with
$\|m\|_\infty\le M$. For $N\ge2$ and any $\theta^N\in\Theta^N$, set
$M_N=N^{-1}\sum_{j=1}^N\mu_{\theta_j}[m]$. Then
\[
\max_{1\le i\le N}\varepsilon_{i,N}(\theta^N;m)
\le C_M\left(\frac{1+\gamma}{N}+\|M_N-m\|_{L^2(0,T)}^2\right).
\]
The constant $C_M$ is uniform over the type vectors.
If $F_N=N^{-1}\sum_{j=1}^N\delta_{\theta_j}$ converges weakly to $F$,
the limit superior of the maximum gain, for fixed $m$ and $\gamma$,
is bounded by $C_M\|\mathscr T[m]-m\|_{L^2(0,T)}^2$.
In particular, at any mean field fixed point $m^*$ with
$\|m^*\|_\infty\le M$, the maximum gain tends to zero. The bound
is $O(N^{-1})$ if $\|M_N-m^*\|_{L^2(0,T)}=O(N^{-1/2})$,
with $M_N$ computed at $m^*$.
\end{corollary}

\begin{corollary}
\label{cor:sampled-nash}
Under the assumptions of Theorem~\ref{thm:main}, let $m^*$ be a
mean field fixed point with $\|m^*\|_\infty\le M$. For each $N\ge2$,
$\gamma>0$ and $0<\delta<1$, with $F^{\otimes N}$-probability at
least $1-\delta$,
\[
\max_{1\le i\le N}\varepsilon_{i,N}(\theta^N;m^*)
\le \frac{C_M}{N}\bigl(1+\gamma+\log(1/\delta)\bigr).
\]
The constant is independent of $N$, $\gamma$ and $\delta$.
\end{corollary}

We express the gain through the candidate and best-response
densities. The candidate's state moments then give the fixed-type
estimate used in the theorem and both corollaries.

\subsection{The Finite-Player Response and Its Value Gap}
\label{subsec:finite-response}

We use the type-$\theta_j$ coefficients from Section~3 and write
$X_j=X^j$, $W_j=W^j$, and $k_j=k_j[m]$. The candidate mean is
$u_j^m(t,x)=-K_j(t)x-k_j$.
When player $i$ uses a density of mean $u_i$ and variance $v_i$,
its state satisfies
\begin{equation}
\label{eq:deviating-state}
dX_i={}(A_iX_i+\barA_i\dmN+B_i u_i)dt +\sqrt{(\sigma_i+C_iX_i+\barC_i\dmN+D_i u_i)^2+D_i^2v_i}\,dW_i.
\end{equation}
For each $j\ne i$, the opponent retains
$\mathcal N(-K_jX_j-k_j,\gamma/\widehat R_j)$, so
\begin{equation}
\label{eq:opponent-state}
\begin{aligned}
dX_j={}&(A_j^{\rm cl}X_j+\barA_j\dmN-B_jk_j)dt\\
&+\sqrt{((C_j-D_jK_j)X_j+\barC_j\dmN+\sigma_j-D_jk_j)^2
+D_j^2\gamma/\widehat R_j}\,dW_j.
\end{aligned}
\end{equation}
Through $\dmN$, player $i$'s deviation can affect the opponents'
states. The full-information control problem has state
$(X_1,\ldots,X_N)$.

Let $V_i^N(t,x)$ denote this unilateral value from
$x=(x_1,\ldots,x_N)$, with discount measured from $t$.
For $j\ne i$, let $a_j^N(t,x)$ be the variance in \eqref{eq:opponent-state}
evaluated at $x$. The verification HJB is
\begin{equation}
\label{eq:component-finite-HJB}
\begin{aligned}[b]
&-\partial_tV_i^N+\beta V_i^N={}
(A_ix_i+\barA_i\bar x_N)\partial_iV_i^N +\sum_{j\ne i}(A_j^{\rm cl}x_j+\barA_j\bar x_N-B_jk_j)\partial_jV_i^N \\
&+\frac12\sum_{j\ne i}a_j^N(t,x)\partial_{jj}V_i^N +\frac{Q_i}2x_i^2+\frac{\bar Q_i}2\bar x_N^2 \\
&+\inf_\pi\biggl\{\int\biggl[
\frac{R_i}2a^2+B_i\partial_iV_i^N a +\frac{\partial_{ii}V_i^N}2
(\sigma_i+C_ix_i+\barC_i\bar x_N+D_i a)^2\biggr]\pi(da)
-\gamma\mathcal H(\pi)\biggr\}, \\
&V_i^N(T,x)={}\frac{H_i}2x_i^2+\frac{\bar H_i}2\bar x_N^2.
\end{aligned}
\end{equation}
The infimum concerns player $i$'s density. The verification covers
all admissible full-information deviations, including history-dependent rules.

We take the symmetric quadratic ansatz
\[
V_i^N(t,x)=\tfrac12x^\top S(t)x+r(t)^\top x+\kappa_i^N(t).
\]
Its derivatives are
$\partial_iV_i^N=(Sx+r)_i$ and $\partial_{jj}V_i^N=S_{jj}$.
Here $S,r$ replace the limiting coefficients $P_i,\eta_i$.
We suppress the indices $i,\theta^N$ and input $m$.

To match the HJB coefficients, we write $\mathcal A x+b$ for the
drift before adding player $i$'s action and
$\ell_j^\top x+\sigma_j^{(i)}$ for the corresponding affine diffusion term.
Let $\barA^N=(\barA_1,\ldots,\barA_N)^\top$ and set
\begin{equation}
\label{eq:finite-coefficients}
\begin{gathered}
\mathcal A_{jk}
=\bigl(A_j-\mathbf1_{\{j\ne i\}}B_jK_j\bigr)\mathbf1_{\{j=k\}}
+\barA_j/N,\qquad
b_j=-\mathbf1_{\{j\ne i\}}B_jk_j,\\
\ell_j=(C_j-\mathbf1_{\{j\ne i\}}D_jK_j)e_j+\barC_j\mathbf1/N,
\qquad \sigma_j^{(i)}=\sigma_j-\mathbf1_{\{j\ne i\}}D_jk_j.
\end{gathered}
\end{equation}
For an opponent, the variance is
$(\ell_j^\top x+\sigma_j^{(i)})^2+D_j^2\gamma/\widehat R_j$. For player $i$, the diffusion
coefficient under a strict action $a$ is $\ell_i^\top x+\sigma_i+D_i a$.
The state-cost Hessians are
\[
\mathcal Q=Q_i e_ie_i^\top+\bar Q_i\mathbf1\mathbf1^\top/N^2,
\qquad
\mathsf H=H_i e_ie_i^\top+\bar H_i\mathbf1\mathbf1^\top/N^2.
\]

Expanding the own-action expression in the HJB gives
\begin{equation}
\label{eq:finite-action-expansion}
\begin{aligned}
&\frac{R_i}2a^2+B_i(Sx+r)_i a
+\frac{S_{ii}}2(\ell_i^\top x+\sigma_i+D_i a)^2\\
&\qquad=\frac{S_{ii}}2(\ell_i^\top x+\sigma_i)^2
+\frac{\widehat R_i^N}2a^2+(\Lambda^\top x+\lambda)a,
\end{aligned}
\end{equation}
where
\begin{equation}
\label{eq:finite-action-coefficients}
\widehat R_i^N=R_i+D_i^2S_{ii},\qquad
\Lambda=B_iSe_i+D_iS_{ii}\ell_i,\qquad
\lambda=B_i r_i+D_iS_{ii}\sigma_i.
\end{equation}
When $\widehat R_i^N>0$, the Gibbs identity minimizes the last two terms
against the density and its entropy. The minimized own-action term is
$S_{ii}(\ell_i^\top x+\sigma_i)^2/2
-(\Lambda^\top x+\lambda)^2/(2\widehat R_i^N)
-\frac\gamma2\log(2\pi\gamma/\widehat R_i^N)$.
Matching the quadratic, linear and constant coefficients in
\eqref{eq:component-finite-HJB} gives
\begin{equation}
\label{eq:finite-Riccati}
\dot S={}\beta S-\mathcal A^\top S-S\mathcal A
-\sum_j S_{jj}\ell_j\ell_j^\top
+\frac{\Lambda\Lambda^\top}{\widehat R_i^N}-\mathcal Q,\qquad S(T)={}\mathsf H,
\end{equation}
\begin{equation}
\label{eq:finite-affine}
\begin{aligned}
\dot r={}&\beta r-\mathcal A^\top r-Sb
-\sum_jS_{jj}\sigma_j^{(i)}\ell_j+\frac{\Lambda\lambda}{\widehat R_i^N},
\qquad r(T)=0,
\end{aligned}
\end{equation}
The constant has terminal value $\kappa_i^N(T)=0$ and satisfies
\begin{equation}
\label{eq:finite-constant}
\dot\kappa_i^N=\beta\kappa_i^N-r^\top b
-\tfrac12\sum_jS_{jj}(\sigma_j^{(i)})^2
-\tfrac\gamma2\sum_{j\ne i}S_{jj}D_j^2/\widehat R_j
+\frac{\lambda^2}{2\widehat R_i^N}
+\frac\gamma2\log\frac{2\pi\gamma}{\widehat R_i^N}.
\end{equation}

We now verify that \eqref{eq:finite-Riccati}--\eqref{eq:finite-constant}
give the full-information best response to the opponents' fixed feedbacks.

\begin{proposition}
\label{prop:finite-response}
Equations \eqref{eq:finite-Riccati}--\eqref{eq:finite-constant}
have a solution on $[0,T]$, with $S$ positive semidefinite and
$\widehat R_i^N\ge\inf_\theta R_\theta$. The full-information best response is
\begin{equation}
\label{eq:finite-response}
\pi_i^{*,N}(da\mid t,x)
=\mathcal N\left(u_i^{*,N}(t,x),\frac\gamma{\widehat R_i^N(t)}\right)(da),
\qquad u_i^{*,N}=-\frac{\Lambda^\top x+\lambda}{\widehat R_i^N}.
\end{equation}
For every admissible deviation,
\begin{equation}
\label{eq:finite-verification}
J_i^N(\pi_i,\widehat\pi_{-i}^N[m];\theta^N)
-V_i^N(0,x_0\mathbf1)
=\gamma\bE_W\int_0^T e^{-\beta t}
D_{\rm KL}(\pi_{i,t}\|\pi_i^{*,N}(\cdot\mid t,X_t))dt.
\end{equation}
Moreover, $S$ is independent of $m,\gamma$, and $r$ is independent
of $\gamma$.
Put $\varrho_i=\widehat R_i^N/\widehat R_i$. The full-information gain
under the candidate profile is
\begin{equation}
\label{eq:exact-gain}
\begin{aligned}
\varepsilon_{i,N}(\theta^N;m)
={}&\frac12\bE_W\int_0^T e^{-\beta t}\widehat R_i^N
|u_i^{*,N}(t,X_t)-u_i^m(t,X_i)|^2dt\\
&+\frac\gamma2\int_0^T e^{-\beta t}
(\varrho_i-1-\log\varrho_i)dt.
\end{aligned}
\end{equation}
\end{proposition}
\begin{proof}
Since $\widehat R_i^N(T)>0$, we solve \eqref{eq:finite-Riccati} locally
backward. To show that $S$ stays positive semidefinite, we test the
backward derivative at $S\ge0$ on a vector $v$ with $Sv=0$.
The drift and discount terms vanish in this quadratic form, and
$\Lambda^\top v=D_iS_{ii}\ell_i^\top v$.
\[
\begin{aligned}
-v^\top\dot S v
&=\sum_jS_{jj}(\ell_j^\top v)^2
-\frac{D_i^2S_{ii}^2}{\widehat R_i^N}(\ell_i^\top v)^2
+v^\top\mathcal Qv\\
&=\sum_{j\ne i}S_{jj}(\ell_j^\top v)^2
+\frac{R_iS_{ii}}{\widehat R_i^N}(\ell_i^\top v)^2
+v^\top\mathcal Qv\ge0.
\end{aligned}
\]
Adding $\epsilon I$ to the backward right-hand side and the terminal
matrix makes the derivative strictly positive in every null direction.
At first contact with the boundary it would be nonpositive. Passing to
the limit by continuous dependence gives $S\ge0$ and $\widehat R_i^N\ge R_i$.

To continue the solution, we bound its trace using
\[
\|\mathcal A\|_{\rm op}
\le\max_j(|A_j|+|B_jK_j|)+|\barA^N|/\sqrt N\le C,
\qquad \max_j|\ell_j|\le C.
\]
Taking the trace in backward time gives
\[
\begin{aligned}
\frac d{d\tau}\tr S(T-\tau)
&=-\beta\tr S+2\tr(S\mathcal A)
+\sum_jS_{jj}|\ell_j|^2-|\Lambda|^2/\widehat R_i^N+\tr\mathcal Q\\
&\le C\tr S+\tr\mathcal Q,
\qquad \tr\mathcal Q=Q_i+\bar Q_i/N\le C.
\end{aligned}
\]
Since $\tr\mathsf H=H_i+\bar H_i/N\le C$, Gronwall bounds
$\tr S$ and hence $\|S\|_{\rm op}$ up to the end of the local interval.
This bound and $\widehat R_i^N\ge R_i$ extend the solution to $[0,T]$.
Once $S$ is fixed, the equation for $r$ is linear because
$\lambda=B_i r_i+D_iS_{ii}\sigma_i$. The constant equation is then
linear as well.

For the verification, we fix an admissible deviation with action
mean $u_i$ and variance $v_i$. The state equation
\eqref{eq:finite-game-state} has drift $\mathcal A X+b+B_i e_i u_i$
and covariance $\operatorname{diag}_j a_{j,t}^N$. Since
$\nabla V_i^N=SX+r$ and $\partial_{jj}V_i^N=S_{jj}$, the HJB and
the Gibbs identity give
\begin{equation}
\label{eq:finite-KL-residual}
\begin{aligned}
&\partial_tV_i^N-\beta V_i^N
+(\mathcal A X+b+B_i e_i u_i)^\top(SX+r)
+\frac12\sum_j a_{j,t}^N S_{jj}\\
&\quad+\frac{Q_i}2X_i(t)^2+\frac{\bar Q_i}2\dmN(t)^2
+\frac{R_i}2\int a^2\pi_{i,t}(da)-\gamma\mathcal H(\pi_{i,t})\\
&=\gamma D_{\rm KL}(\pi_{i,t}\|\pi_i^{*,N}(\cdot\mid t,X_t)).
\end{aligned}
\end{equation}
Here $a_{i,t}^N=(\sigma_i+C_iX_i+\barC_i\dmN+D_i u_i)^2+D_i^2v_i$;
for $j\ne i$ it is the variance in \eqref{eq:opponent-state}.
The inequality $\dmN^2\le |X|^2/N$ and boundedness of
$K,k,\widehat R^{-1}$ imply
\[
|\mathcal A X+b+B_i e_i u_i|^2+\sum_j a_{j,t}^N
\le C_{N,M,\gamma}\left(1+|X_t|^2+\int a^2\pi_{i,t}(da)\right).
\]
Before taking expectations in the discounted It\^o identity, we
use the stops \eqref{eq:scalar-stopping} with the vector norm and
$\pi_s=\pi_{i,s}$. Nonexplosion and the action-moment condition give
$\tau_n\uparrow T$. In the moment estimate leading to
\eqref{eq:verification-moment}, the squared drift is now
$|\mathcal A X+b+B_i e_i u_i|^2$ and the quadratic variation density
is $\sum_j a_{j,t}^N$. The growth bound gives the same stopped
integral inequality. Gronwall and Fatou yield
\begin{equation}
\label{eq:finite-verification-moment}
\bE_W\sup_{t\le T}|X_t|^2
\le C_{N,M,\gamma}\left(1+|X_0|^2+
\bE_W\int_0^T\int a^2\pi_{i,t}(da)dt\right)<\infty.
\end{equation}
The stochastic integral in the discounted It\^o identity is
\[
\mathcal M_s=\sum_j\int_0^s e^{-\beta t}(SX_t+r)_j
\sqrt{a_{j,t}^N}\,dW_j.
\]
The growth bound and \eqref{eq:finite-verification-moment} imply
$\bE_W\int_0^T\sum_j a_{j,t}^Ndt<\infty$.
By BDG and Cauchy--Schwarz as in \eqref{eq:verification-martingale-bound},
\[
\bE_W\sup_{s\le T}|\mathcal M_s|
\le C_{N,M}\left[\bE_W\left(1+\sup_{t\le T}|X_t|^2\right)\right]^{1/2}
\left[\bE_W\int_0^T\sum_j a_{j,t}^Ndt\right]^{1/2}<\infty.
\]
Thus $\mathcal M$ is a uniformly integrable martingale and
$\bE_W\mathcal M_T=0$.
It\^o's formula and \eqref{eq:finite-KL-residual} give
\[
\begin{aligned}
&e^{-\beta T}\left[\frac{H_i}2X_i(T)^2+
\frac{\bar H_i}2\dmN(T)^2\right]-V_i^N(0,x_0\mathbf1)\\
&+\int_0^T e^{-\beta t}
\left[\frac{Q_i}2X_i(t)^2+\frac{\bar Q_i}2\dmN(t)^2
+\frac{R_i}2\int a^2\pi_{i,t}(da)-\gamma\mathcal H(\pi_{i,t})\right]dt\\
&=\gamma\int_0^T e^{-\beta t}
D_{\rm KL}(\pi_{i,t}\|\pi_i^{*,N}(\cdot\mid t,X_t))dt+\mathcal M_T.
\end{aligned}
\]
The endpoint is integrable by \eqref{eq:finite-verification-moment};
the absolute running cost is bounded by
$C(|X_t|^2+\int a^2\pi_{i,t}(da))+\gamma|\mathcal H(\pi_{i,t})|$.
Taking expectations therefore proves \eqref{eq:finite-verification}.
Under \eqref{eq:finite-response} the drift is affine, and the
reverse triangle inequality bounds the diffusion Lipschitz
constants by $|\ell_i-D_i\Lambda/\widehat R_i^N|$ for player $i$
and $|\ell_j|$ for each opponent. The joint SDE has a unique strong
solution with finite maximal second moment. The response has
action moment $|u_i^{*,N}(t,X_t)|^2+\gamma/\widehat R_i^N
\le C_{N,m}(1+|X_t|^2)+C\gamma$.
The entropy $\frac12\log(2\pi e\gamma/\widehat R_i^N)$ is bounded
in time. The response is admissible, and its KL term is zero.

The verification also holds from any starting time and deterministic
state vector: the same moment and martingale estimates apply with
that initial condition. Finally, the coefficients of
\eqref{eq:finite-Riccati} do not use $m$ or $\gamma$; in
\eqref{eq:finite-affine} only $m$ occurs through $k_j$.
The ODE solutions depend continuously on the type coefficients.

Under the candidate profile, all action means satisfy
\eqref{eq:uniform-affine-coefficients}. Write $\widehat\pi_{i,t}^N$
for the candidate density along $X^i$. The moment calculation gives
\[
\sup_{\theta^N}\bE_W\sup_{t\le T}|X_t|^2\le C_{N,M}(1+\gamma),
\qquad
\bE_i^\zeta\int_0^T\int a^2\widehat\pi_{i,t}^N(da)dt
\le C_{N,M}(1+\gamma).
\]
Its Gaussian entropy is bounded uniformly over types for fixed
$\gamma$. Thus the candidate feedback is admissible and
has finite type-averaged cost, so $\mathcal E_{i,N}(\zeta;m)\ge0$.

For Gaussian laws with means $u_m,u_*$ and variances
$v_m,v_*>0$, direct integration of their log densities gives
\[
D_{\rm KL}(\mathcal N(u_m,v_m)\|\mathcal N(u_*,v_*))
=\tfrac12\left[\frac{(u_m-u_*)^2}{v_*}
+\frac{v_m}{v_*}-1-\log\frac{v_m}{v_*}\right].
\]
Taking $v_m=\gamma/\widehat R_i$, $v_*=\gamma/\widehat R_i^N$ in
\eqref{eq:finite-verification}, under the candidate profile, gives \eqref{eq:exact-gain}.
\end{proof}

These moment bounds justify verification at fixed $N$ and $\gamma$.
The next subsection estimates the gain uniformly in $N$.
% performance.tex
\subsection{Mean and Covariance Comparison}
\label{subsec:mean-covariance-comparison}

To estimate the mean-feedback term in \eqref{eq:exact-gain},
we write the feedback difference in affine form.
\begin{equation}
\label{eq:feedback-difference-coefficients}
u_i^{*,N}(t,x)-u_i^m(t,x_i)
=\omega_i^\top x-\lambda/\widehat R_i^N+k_i,\qquad
\omega_i=-\Lambda/\widehat R_i^N+K_i e_i.
\end{equation}
The affine difference is independent of $\gamma$, and
$\omega_i$ is independent of $m$, by Proposition~\ref{prop:finite-response}.
Let $X$ be the finite-game state under $\widehat\pi^N[m]$, and let
$\nu=\bE_W X$ and $\Sigma_\gamma=\operatorname{Cov}_W(X)$ be its mean
and covariance at fixed types. Centering at $\nu$ gives
\begin{equation}
\label{eq:feedback-second-moment}
\bE_W|u_i^{*,N}(t,X)-u_i^m(t,X_i)|^2
=\omega_i^\top\Sigma_\gamma\omega_i
+|u_i^{*,N}(t,\nu)-u_i^m(t,\nu_i)|^2.
\end{equation}
The covariance term depends on the slope $\omega_i$. Bounds on the
entries of $S$ control both this slope and the variance ratio in
\eqref{eq:exact-gain}. For the second term, we compare the feedbacks
at the type means in \eqref{eq:mu-theta} before returning to $\nu$.

\begin{lemma}
\label{lem:row-bounds}
With $C$ independent of $N,i,\theta^N,m,\gamma$, the response matrix
satisfies, for $j,k\ne i$,
\begin{equation}
\label{eq:row-bounds}
0\le S_{ii}\le C,\quad
S_{jj}\le C/N^2,\quad
|S_{ij}|\le C/N,\quad
|S_{jk}|\le C/N^2,
\end{equation}
and
\begin{equation}
\label{eq:diagonal-comparison}
\sup_{t\le T}|S_{ii}(t)-P_i(t)|\le C/N.
\end{equation}
\end{lemma}
\begin{proof}
We first bound the opponents' diagonal entries. The cross entries
then follow from $S\ge0$. The proof of Proposition~\ref{prop:finite-response}
gives $\sup_{t\le T}\tr S(t)\le C$. For $j\ne i$, \eqref{eq:finite-Riccati} gives
\[
\begin{aligned}
-\dot S_{jj}
&=(2A_j^{\rm cl}-\beta)S_{jj}
+\frac2N\sum_l\barA_l S_{lj}\\
&\quad+\sum_l S_{ll}(\ell_l)_j^2
-\frac{\Lambda_j^2}{\widehat R_i^N}+\frac{\bar Q_i}{N^2},
\qquad S_{jj}(T)=\frac{\bar H_i}{N^2}.
\end{aligned}
\]
Positive semidefiniteness implies
$|S_{ij}|\le C\sqrt{S_{jj}}$ and
$|S_{lj}|\le\max_{k\ne i}S_{kk}$ for $l,j\ne i$.
Consequently, Young's inequality gives
\[
\begin{aligned}
\frac2N\left|\sum_l\barA_l S_{lj}\right|
&\le\frac2N|\barA_i S_{ij}|+\frac2N\sum_{l\ne i}|\barA_l S_{lj}|\\
&\le \frac C N\sqrt{S_{jj}}+C\max_{k\ne i}S_{kk}
\le \frac C{N^2}+C\max_{k\ne i}S_{kk}.
\end{aligned}
\]
The explicit components of $\ell_l$ give
\[
\begin{aligned}
\sum_lS_{ll}(\ell_l)_j^2
&=S_{jj}(C_j-D_jK_j+\barC_j/N)^2
+\frac1{N^2}\sum_{l\ne j}S_{ll}\barC_l^2\\
&\le 2(C_j-D_jK_j)^2S_{jj}+\frac2{N^2}\sum_lS_{ll}\barC_l^2
\le C S_{jj}+\frac C{N^2}\tr S
\le C\max_{k\ne i}S_{kk}+\frac C{N^2}.
\end{aligned}
\]
The term $-\Lambda_j^2/\widehat R_i^N$ is nonpositive.
Integrating the resulting upper bound from $t$ to $T$ and taking
the maximum over $j\ne i$ yields
\[
\max_{j\ne i}S_{jj}(t)
\le \frac C{N^2}
+C\int_t^T\max_{j\ne i}S_{jj}(s)ds.
\]
Backward Gronwall gives $S_{jj}\le C/N^2$ for every $j\ne i$.
Since $S_{ii}\le\tr S\le C$, positive semidefiniteness also gives
\[
|S_{ij}|^2\le S_{ii}S_{jj}\le C/N^2,\qquad
|S_{jk}|^2\le S_{jj}S_{kk}\le C/N^4\quad(j,k\ne i),
\]
which proves \eqref{eq:row-bounds}.

For the own-state curvature, we subtract \eqref{eq:asy-Riccati}
from the $(i,i)$ entry of \eqref{eq:finite-Riccati}.
\begin{equation}
\label{eq:diagonal-difference}
\begin{aligned}
\frac d{dt}(S_{ii}-P_i)
={}&(\beta-2A_i-C_i^2)(S_{ii}-P_i)\\
&+(B_i+D_iC_i)^2\left(\frac{S_{ii}^2}{R_i+D_i^2S_{ii}}
-\frac{P_i^2}{R_i+D_i^2P_i}\right)+\mathcal R_N,
\end{aligned}
\end{equation}
with $(S_{ii}-P_i)(T)=\bar H_i/N^2$ and
\[
\begin{aligned}
\mathcal R_N={}&-\frac2N e_i^\top S\barA^N
-S_{ii}\left[(C_i+\barC_i/N)^2-C_i^2\right]
-\frac{\sum_{j\ne i}S_{jj}\barC_j^2+\bar Q_i}{N^2}\\
&+\frac{S_{ii}^2}{\widehat R_i^N}
\left[(B_i+D_i(C_i+\barC_i/N))^2-(B_i+D_iC_i)^2\right].
\end{aligned}
\]
The row bounds give $|e_i^\top S\barA^N|\le C$ and
$\sum_jS_{jj}\barC_j^2\le C$. Both differences of squares
are bounded by $C/N$, so $|\mathcal R_N|\le C/N$.
On the bounded nonnegative range of $S_{ii}$ and $P_i$,
\[
\frac d{ds}\frac{s^2}{R_i+D_i^2s}
=\frac{s(2R_i+D_i^2s)}{(R_i+D_i^2s)^2}\le C.
\]
The integral form of the difference equation then gives
\[
|S_{ii}(t)-P_i(t)|
\le\frac{\bar H_i}{N^2}
+C\int_t^T\left(|S_{ii}(s)-P_i(s)|+N^{-1}\right)ds.
\]
Backward Gronwall proves \eqref{eq:diagonal-comparison}.
\end{proof}

Equations \eqref{eq:finite-action-coefficients} and
\eqref{eq:feedback-difference-coefficients} give the entries of $\omega_i$.
\[
(\omega_i)_i=-\frac{(B_i+D_iC_i)R_i(S_{ii}-P_i)}{\widehat R_i\widehat R_i^N}
-\frac{D_iS_{ii}\barC_i}{N\widehat R_i^N},\;
(\omega_i)_j=-\frac{B_iS_{ij}+D_iS_{ii}\barC_i/N}{\widehat R_i^N},\quad j\ne i.
\]
The first identity uses
$P_i/\widehat R_i-S_{ii}/\widehat R_i^N
=-R_i(S_{ii}-P_i)/(\widehat R_i\widehat R_i^N)$.
Lemma~\ref{lem:row-bounds} therefore gives
$\max_j|(\omega_i)_j|\le C/N$, and hence
\begin{equation}
\label{eq:feedback-coefficient-bounds}
|\omega_i|^2\le C/N^2+C(N-1)/N^2\le C/N,\qquad
\sum_j|(\omega_i)_j|\le C/N+C(N-1)/N\le C.
\end{equation}

We now prove Theorem~\ref{thm:main}.
\begin{proof}
The action weights are bounded above and away from zero.
Equation \eqref{eq:diagonal-comparison} therefore gives
\[
|\varrho_i-1|=\frac{D_i^2|S_{ii}-P_i|}{\widehat R_i}\le C/N.
\]
Taylor's formula at $1$ gives
\begin{equation}
\label{eq:variance-ratio-bound}
0\le\varrho_i-1-\log\varrho_i
\le\frac{|\varrho_i-1|^2}{2\min(1,\varrho_i)^2}\le C/N^2,
\end{equation}
since the ratio stays in a compact subset of $(0,\infty)$.
The variance-ratio contribution to \eqref{eq:exact-gain} is
at most $C\gamma/N^2$. We now estimate the covariance term
of the remaining mean-feedback difference in
\eqref{eq:feedback-second-moment}. Set
\[
\mathcal A^{\rm cl}=\operatorname{diag}(A_j^{\rm cl})
+\barA^N\mathbf1^\top/N,\qquad
\ell_j^{\rm cl}=(C_j-D_jK_j)e_j+\barC_j\mathbf1/N.
\]
Thus $\mathcal A^{\rm cl}=\mathcal A-B_iK_i e_ie_i^\top$,
$\ell_j^{\rm cl}=\ell_j$ for $j\ne i$, and
$\ell_i^{\rm cl}=\ell_i-D_iK_i e_i$.
The coefficient definitions give
\begin{equation}
\label{eq:implemented-coefficient-bounds}
\begin{gathered}
\|\mathcal A^{\rm cl}\|_\infty\le\max_j(|A_j^{\rm cl}|+|\barA_j|)\le C,
\quad\|\mathcal A^{\rm cl}\|_{\rm op}
\le\max_j|A_j^{\rm cl}|+|\barA^N|/\sqrt N\le C,\\
\max_j\sum_k|(\ell_j^{\rm cl})_k|
\le\max_j(|C_j-D_jK_j|+|\barC_j|)\le C.
\end{gathered}
\end{equation}
In particular $\max_j|\ell_j^{\rm cl}|\le C$.
The moment bound \eqref{eq:finite-verification-moment} makes
the state stochastic integrals square integrable. Taking
expectations under the candidate profile gives
\begin{equation}
\label{eq:finite-mean-ODE}
\dot\nu=\mathcal A^{\rm cl}\nu-(B_jk_j)_j,\qquad
\nu(0)=x_0\mathbf1.
\end{equation}
By \eqref{eq:uniform-affine-coefficients},
\[
\|\nu(t)\|_\infty\le |x_0|+C_M t
+C\int_0^t\|\nu(s)\|_\infty ds.
\]
Gronwall gives $\sup_{t\le T}\|\nu(t)\|_\infty\le C_M$,
independently of $\gamma$.

For the covariance calculation, we write $Z=X-\nu$ and denote the
affine diffusion term at the mean by $f_j=(\ell_j^{\rm cl})^\top\nu+\sigma_j-D_jk_j$.
The mean bound and \eqref{eq:implemented-coefficient-bounds} give
\[
|f_j|\le\sum_k|(\ell_j^{\rm cl})_k|\,\|\nu\|_\infty
+|\sigma_j|+|D_jk_j|\le C_M.
\]
In the original realization \eqref{eq:finite-game-state},
the centered coordinates satisfy
\[
dZ_j=(\mathcal A^{\rm cl}Z)_jdt+\sqrt{a_{j,t}^N}\,dW_j,\quad Z(0)=0,\;
a_{j,t}^N=((\ell_j^{\rm cl})^\top Z+f_j)^2+\gamma D_j^2/\widehat R_j.
\]
For the independent Brownian coordinates, It\^o's product formula gives
\[
\begin{aligned}
d(Z_jZ_k)
&=\left[(\mathcal A^{\rm cl}Z)_jZ_k+
Z_j(\mathcal A^{\rm cl}Z)_k+\mathbf1_{\{j=k\}}a_{j,t}^N\right]dt\\
&\quad+Z_k\sqrt{a_{j,t}^N}\,dW_j+Z_j\sqrt{a_{k,t}^N}\,dW_k.
\end{aligned}
\]
The moment bound gives $\bE_W\sup_{t\le T}|Z_t|^2<\infty$ and $\bE_W\int_0^T\sum_j a_{j,t}^Ndt<\infty$.
By BDG, the stochastic integrals in the product formula have
integrable maxima, since
\[
\bE_W\left(\int_0^T|Z_t|^2\sum_j a_{j,t}^Ndt\right)^{1/2}
\le\left[\bE_W\sup_{t\le T}|Z_t|^2\right]^{1/2}
\left[\bE_W\int_0^T\sum_j a_{j,t}^Ndt\right]^{1/2}<\infty.
\]
They therefore have zero expectation. Since $\bE_W Z=0$,
\[
\bE_W((\ell_j^{\rm cl})^\top Z+f_j)^2
=(\ell_j^{\rm cl})^\top\Sigma_\gamma\ell_j^{\rm cl}+f_j^2.
\]
Taking expectations in the product formula yields the centered second-moment
equation (cf. \cite[Equation~(2.11)]{GRZ24})
\begin{equation}
\label{eq:covariance-ODE}
\dot\Sigma_\gamma=\mathcal A^{\rm cl}\Sigma_\gamma
+\Sigma_\gamma(\mathcal A^{\rm cl})^\top
+\operatorname{diag}_j\left\{(\ell_j^{\rm cl})^\top\Sigma_\gamma\ell_j^{\rm cl}
+f_j^2+\gamma D_j^2/\widehat R_j\right\},\qquad
\Sigma_\gamma(0)=0.
\end{equation}

Since $\Sigma_\gamma$ is positive semidefinite, for every unit
vector $v$ its integral equation gives
\[
\begin{aligned}
v^\top\Sigma_\gamma(t)v
&\le\int_0^t\Bigl[
\bigl(2\|\mathcal A^{\rm cl}\|_{\rm op}
+\max_j|\ell_j^{\rm cl}|^2\bigr)\|\Sigma_\gamma\|_{\rm op}
+\max_j f_j^2+\gamma\max_jD_j^2/\widehat R_j\Bigr]ds\\
&\le C\int_0^t\|\Sigma_\gamma(s)\|_{\rm op}ds+C_M(1+\gamma)t.
\end{aligned}
\]
Taking the supremum over unit vectors and applying Gronwall gives
\begin{equation}
\label{eq:covariance-bound}
\sup_{t\le T}\|\Sigma_\gamma(t)\|_{\rm op}\le C_M(1+\gamma).
\end{equation}
The constants in this estimate are independent of $N$.
By \eqref{eq:feedback-coefficient-bounds},
\begin{equation}
\label{eq:covariance-contribution}
\int_0^T e^{-\beta t}\widehat R_i^N
\omega_i^\top\Sigma_\gamma\omega_i\,dt
\le C\int_0^T\|\Sigma_\gamma\|_{\rm op}|\omega_i|^2dt
\le C_M(1+\gamma)/N.
\end{equation}

To estimate the second term in \eqref{eq:feedback-second-moment}, we
compare $\nu$ with the type means $\mu=(\mu_1,\ldots,\mu_N)^\top$,
where $\mu_j=\mu_{\theta_j}[m]$ solves \eqref{eq:mu-theta}.
Their average and its discrepancy from $m$ are
\begin{equation}
\label{eq:empirical-mean-residual}
M_N=N^{-1}\sum_j\mu_j,\qquad d_N=M_N-m.
\end{equation}
The residual $d_N$ records both mean inconsistency and type composition. Subtracting
\eqref{eq:mu-theta} from \eqref{eq:finite-mean-ODE} gives
\[
\frac d{dt}(\nu-\mu)=\mathcal A^{\rm cl}(\nu-\mu)+\barA^N d_N,
\qquad (\nu-\mu)(0)=0.
\]
The integral equation and $\|\mathcal A^{\rm cl}\|_\infty\le C$
give
$\|\nu(t)-\mu(t)\|_\infty
\le C\int_0^t(\|\nu(s)-\mu(s)\|_\infty+|d_N(s)|)ds$.
Gronwall and Cauchy--Schwarz yield
\begin{equation}
\label{eq:finite-mean-comparison}
\sup_{t\le T}\|\nu(t)-\mu(t)\|_\infty^2\le C \|d_N\|_{L^2(0,T)}^2.
\end{equation}
At $\mu$, we compare the feedbacks using the value gradient difference
\begin{equation}
\label{eq:centered-affine}
h_i=P_i\mu_i+\eta_i,\qquad z=S\mu+r-h_i e_i.
\end{equation}
Here $h_i=\partial_x V^{\theta_i}(t,\mu_i;m)$, and
\[
z_i=\partial_iV_i^N(t,\mu)-h_i,\qquad
z_j=\partial_jV_i^N(t,\mu)\quad(j\ne i).
\]
Let $\sigma_j^\mu=\sigma_j+C_j\mu_j+\barC_j m+D_j u_j^m(t,\mu_j)$.
To obtain the equation for $h_i$, we differentiate the scalar HJB
at $x=\mu_i$. The optimal variance is constant in $x$, and the
Hamiltonian is stationary in the action mean at the optimum.
\[
-\partial_t V_x^{\theta_i}(t,\mu_i;m)+\beta h_i
=A_i h_i+P_i\dot\mu_i+Q_i\mu_i+C_iP_i\sigma_i^\mu.
\]
The chain rule $\dot h_i=\partial_t V_x^{\theta_i}(t,\mu_i;m)+P_i\dot\mu_i$ gives
\begin{equation}
\label{eq:mean-costate}
\dot h_i=(\beta-A_i)h_i-Q_i\mu_i-C_iP_i\sigma_i^\mu,
\qquad h_i(T)=H_i\mu_i(T).
\end{equation}
The feedback formulas at these means give
\[
\begin{aligned}
\widehat R_i^N u_i^{*,N}(t,\mu)
&=-B_i(S\mu+r)_i-D_iS_{ii}(\sigma_i+C_i\mu_i+\barC_i M_N),\\
\widehat R_i u_i^m(t,\mu_i)
&=-B_i h_i-D_iP_i(\sigma_i+C_i\mu_i+\barC_i m).
\end{aligned}
\]
We multiply the feedback difference by $\widehat R_i^N$ and use
$\widehat R_i^N-\widehat R_i=D_i^2(S_{ii}-P_i)$ to obtain
\begin{equation}
\label{eq:mean-response-at-mu}
\widehat R_i^N[u_i^{*,N}(t,\mu)-u_i^m(t,\mu_i)]
=-B_i z_i-D_i(S_{ii}-P_i)\sigma_i^\mu-D_iS_{ii}\barC_i d_N.
\end{equation}

To estimate $z$, we differentiate both gradients along $\mu$,
separating the forcing proportional to $d_N$ from the $O(N^{-1})$ terms.
The type-mean equations give
\[
\dot\mu=\mathcal A\mu+b+B_i e_i u_i^m(t,\mu_i)+\barA^N(m-M_N).
\]
We differentiate $S\mu+r$ using
\eqref{eq:finite-Riccati}--\eqref{eq:finite-affine} and the mean equation.
The terms $S\mathcal A\mu$ and $Sb$ cancel, and
$\Lambda(\Lambda^\top\mu+\lambda)/\widehat R_i^N
=-\Lambda u_i^{*,N}(t,\mu)$. Hence
\begin{equation}
\label{eq:mean-gradient-product}
\begin{aligned}
\frac d{dt}(S\mu+r)={}&(\beta I-\mathcal A^\top)(S\mu+r)-\mathcal Q\mu
-\sum_jS_{jj}\ell_j(\ell_j^\top\mu+\sigma_j^{(i)})
-\Lambda u_i^{*,N}(t,\mu)\\
&+B_iSe_i u_i^m(t,\mu_i)-S\barA^N d_N.
\end{aligned}
\end{equation}
For $j\ne i$ we have
$\ell_j^\top\mu+\sigma_j^{(i)}=\sigma_j^\mu+\barC_jd_N$.
For player $i$, the own action has not been inserted in $\ell_i$,
so
\[
\ell_i^\top\mu+\sigma_i^{(i)}
=\sigma_i^\mu-D_i u_i^m(t,\mu_i)+\barC_i d_N.
\]
Consequently,
\[
\sum_jS_{jj}\ell_j(\ell_j^\top\mu+\sigma_j^{(i)})
=\sum_jS_{jj}\ell_j\sigma_j^\mu
+d_N\sum_jS_{jj}\barC_j\ell_j
-D_iS_{ii}\ell_i u_i^m(t,\mu_i).
\]
The last term combines with $B_iSe_i u_i^m$ in
\eqref{eq:mean-gradient-product} to give $\Lambda u_i^m$.
Subtracting the scalar gradient equation removes the own-state drift
and cost terms.
\[
\begin{aligned}
\mathcal A^\top e_i&=A_i e_i+\barA_i\mathbf1/N,\qquad \mathcal Q\mu=Q_i\mu_i e_i+\bar Q_iM_N\mathbf1/N, \\
(\beta I-\mathcal A^\top)h_i e_i
-\dot h_i e_i-\mathcal Q\mu
&=C_iP_i\sigma_i^\mu e_i-\frac{\barA_i h_i+\bar Q_i M_N}{N}\mathbf1.
\end{aligned}
\]
Using \eqref{eq:mean-response-at-mu} in the resulting equation gives
\begin{equation}
\label{eq:z-equation}
\begin{aligned}
\dot z&=\left(\beta I-\mathcal A^\top+
\frac{B_i}{\widehat R_i^N}\Lambda e_i^\top\right)z+F_0+F_1d_N,
&z(T)&=\frac{\bar H_iM_N(T)}N\mathbf1,\\
F_0&=\frac{D_i(S_{ii}-P_i)}{\widehat R_i^N}\Lambda\sigma_i^\mu
-\sum_j S_{jj}\ell_j\sigma_j^\mu+C_iP_i\sigma_i^\mu e_i
-\frac{\barA_i h_i+\bar Q_i M_N}{N}\mathbf1,\\
F_1&=\frac{D_iS_{ii}\barC_i}{\widehat R_i^N}\Lambda
-\sum_j S_{jj}\barC_j\ell_j-S\barA^N.
\end{aligned}
\end{equation}
At $T$, $r(T)=0$ and $h_i(T)=H_i\mu_i(T)$, so
\[
z(T)=\left(H_i e_ie_i^\top+\frac{\bar H_i}{N^2}\mathbf1\mathbf1^\top\right)\mu(T)
-H_i\mu_i(T)e_i=\frac{\bar H_iM_N(T)}N\mathbf1.
\]

The bounds \eqref{eq:uniform-affine-coefficients} make
$\mu_j,h_i,\sigma_j^\mu$ uniformly bounded by $C_M$.
By \eqref{eq:finite-action-coefficients} and Lemma~\ref{lem:row-bounds},
\[
|\Lambda_i|\le C,\qquad |\Lambda_j|\le C/N\ (j\ne i),
\qquad \sum_j|\Lambda_j|\le C.
\]
For any coordinate $k$, the diffusion terms in $F_0$ satisfy
\[
\begin{aligned}
\left(\sum_jS_{jj}\ell_j\sigma_j^\mu-C_iP_i\sigma_i^\mu e_i\right)_k
={}&\mathbf1_{\{k=i\}}C_i(S_{ii}-P_i)\sigma_i^\mu\\
&+\mathbf1_{\{k\ne i\}}S_{kk}(C_k-D_kK_k)\sigma_k^\mu
+\frac1N\sum_jS_{jj}\barC_j\sigma_j^\mu .
\end{aligned}
\]
The first term is bounded by $C_M/N$ by \eqref{eq:diagonal-comparison},
the second by $C_M/N^2$, and the last by $C_M\tr S/N$. For the terms
multiplying $d_N$, the same component formula for $\ell_j$ gives
\[
\left\|\sum_j S_{jj}\barC_j\ell_j\right\|_\infty
\le C\max_j S_{jj}+C\tr S/N\le C,\qquad
\|S\barA^N\|_\infty\le\|\barA\|_\infty\max_k\sum_j|S_{kj}|\le C.
\]
Together with $|S_{ii}-P_i|\le C/N$ and the bounds on $\Lambda$,
these estimates give
\[
\|F_0\|_\infty\le C_M\left(|S_{ii}-P_i|+N^{-1}\right)\le C_M/N,
\qquad \|F_1\|_\infty\le C_M.
\]
Moreover,
$\|\mathcal A^\top\|_\infty
\le\max_j(|A_j|+|B_jK_j|)+N^{-1}\sum_j|\barA_j|\le C$
and $\|\Lambda e_i^\top\|_\infty=\max_j|\Lambda_j|\le C$.
Thus the coefficient matrix in \eqref{eq:z-equation} has induced
infinity norm at most $C$, independently of $N$. Its backward
integral equation gives
\[
\|z(t)\|_\infty\le\frac{C_M}{N}
+C\int_t^T\|z(s)\|_\infty ds
+C_M\int_t^T|d_N(s)|ds.
\]
Backward Gronwall and Cauchy--Schwarz imply
\begin{equation}
\label{eq:z-bound}
\sup_{t\le T}\|z(t)\|_\infty^2\le C_M(N^{-2}+\|d_N\|_{L^2(0,T)}^2).
\end{equation}

The boundedness of $\sigma_i^\mu$ and the positive lower bound
on $\widehat R_i^N$ in \eqref{eq:mean-response-at-mu} imply
\[
|u_i^{*,N}(t,\mu)-u_i^m(t,\mu_i)|^2
\le C_M\bigl(|z_i|^2+|S_{ii}-P_i|^2+|d_N|^2\bigr).
\]
By \eqref{eq:z-bound} and \eqref{eq:diagonal-comparison},
\[
\int_0^T|u_i^{*,N}(t,\mu)-u_i^m(t,\mu_i)|^2dt
\le C_M(N^{-2}+\|d_N\|_{L^2(0,T)}^2).
\]
The affine formula \eqref{eq:feedback-difference-coefficients}
transfers the comparison from $\mu$ to the actual mean $\nu$.
\[
\begin{aligned}
u_i^{*,N}(t,\nu)-u_i^m(t,\nu_i)
={}&u_i^{*,N}(t,\mu)-u_i^m(t,\mu_i)+\omega_i^\top(\nu-\mu).
\end{aligned}
\]
The inequality $|\omega_i^\top(\nu-\mu)|
\le\sum_j|(\omega_i)_j|\,\|\nu-\mu\|_\infty$ gives
\begin{equation}
\label{eq:deterministic-feedback-bound}
\begin{aligned}
\int_0^T|u_i^{*,N}(t,\nu)-u_i^m(t,\nu_i)|^2dt
&\le 2\int_0^T|u_i^{*,N}(t,\mu)-u_i^m(t,\mu_i)|^2dt\\
&\quad+2T\sup_t\left(\sum_j|(\omega_i)_j|\right)^2
\sup_t\|\nu-\mu\|_\infty^2\\
&\le C_M(N^{-2}+\|d_N\|_{L^2(0,T)}^2).
\end{aligned}
\end{equation}
Here we used \eqref{eq:feedback-coefficient-bounds} and
\eqref{eq:finite-mean-comparison}.

We insert \eqref{eq:covariance-contribution} and
\eqref{eq:deterministic-feedback-bound} in
\eqref{eq:feedback-second-moment}. The exact gap
\eqref{eq:exact-gain}, the variance bound
\eqref{eq:variance-ratio-bound} and $\widehat R_i^N\le C$ give
\begin{equation}
\label{eq:fixed-population-bound}
0\le\varepsilon_{i,N}(\theta^N;m)
\le C_M\left(\frac{1+\gamma}{N}+\|d_N\|_{L^2(0,T)}^2\right).
\end{equation}
The ODEs \eqref{eq:finite-Riccati}--\eqref{eq:finite-affine},
\eqref{eq:finite-mean-ODE} and \eqref{eq:covariance-ODE} give
continuous dependence of the feedback coefficients and state moments on $\theta^N$.
Equations \eqref{eq:exact-gain} and \eqref{eq:feedback-second-moment}
then show that the gain is Borel in $\theta^N$.
We now fix $\theta_i=\zeta$ and average over the opponent types.
Their independence gives, at each time,
\begin{equation}
\label{eq:type-variance-identity}
\bE_i^\zeta|M_N-m|^2
={}\frac{N-1}{N^2}\operatorname{Var}_F(\mu_\theta[m]) +\left|\mathscr T[m]-m+
\frac{\mu_\zeta[m]-\mathscr T[m]}N\right|^2.
\end{equation}
The variance term comes from sampling opponent types; the square
contains the consistency residual and the fixed own-type correction.
By \eqref{eq:uniform-affine-coefficients},
$\|\mu_\zeta[m]-\mathscr T[m]\|_\infty\le2C_M$ uniformly over $\zeta$,
independently of $\gamma$. Integrating \eqref{eq:type-variance-identity} gives
\[
\bE_i^\zeta \|d_N\|_{L^2(0,T)}^2
\le C_M/N+2\|\mathscr T[m]-m\|_{L^2(0,T)}^2.
\]
Combining this inequality with the type average of
\eqref{eq:fixed-population-bound} gives \eqref{eq:main-residual}.
Equation~\eqref{eq:own-observation-bound} gives the stated consequence
for decentralized deviations.
\end{proof}

\begin{proof}[Proof of Corollary~\ref{cor:fixed-types}]
Taking the maximum in \eqref{eq:fixed-population-bound} gives the
first bound. For each $t$, weak convergence gives
$M_N(t)\to\mathscr T[m](t)$ because $\theta\mapsto\mu_\theta[m](t)$
is continuous. The uniform bound \eqref{eq:uniform-affine-coefficients}
and dominated convergence give $M_N\to\mathscr T[m]$ in $L^2(0,T)$.
Thus $\|M_N-m\|_{L^2}^2\to\|\mathscr T[m]-m\|_{L^2}^2$, which is
zero at $m=m^*$. The rate follows by squaring the assumed
$O(N^{-1/2})$ discrepancy.
\end{proof}

\begin{proof}[Proof of Corollary~\ref{cor:sampled-nash}]
We take $m=m^*$ in \eqref{eq:empirical-mean-residual} and use
$\bP$ and $\bE$ under the type law $F^{\otimes N}$. Let $B_*>0$ satisfy
$\sup_\theta\|\mu_\theta[m^*]\|_{L^2}\le B_*$ by
\eqref{eq:uniform-affine-coefficients}.
The summands $\mu_{\theta_j}[m^*]-m^*$ are independent and centered
in $L^2(0,T)$. Hence
\[
\bE\|d_N\|_{L^2}^2
=\frac1N\left(\int_\Theta\|\mu_\theta[m^*]\|_{L^2}^2F(d\theta)
-\|m^*\|_{L^2}^2\right)\le\frac{B_*^2}{N}.
\]
In particular, $\bE\|d_N\|_{L^2}\le B_*/\sqrt N$.
By the reverse triangle inequality, replacing $\theta_j$ by $\theta_j'$
changes $\|d_N\|_{L^2}$ by at most
$\|\mu_{\theta_j}[m^*]-\mu_{\theta_j'}[m^*]\|_{L^2}/N$,
hence by at most $2B_*/N$.
For $r>0$, McDiarmid's inequality
\cite[Theorem~3.1]{McDiarmid98} gives
\[
\bP\bigl(\|d_N\|_{L^2}-\bE\|d_N\|_{L^2}\ge r\bigr)
\le \exp\!\left(-\frac{2r^2}{N(2B_*/N)^2}\right)
=\exp\!\left(-\frac{Nr^2}{2B_*^2}\right).
\]
With $r=B_*\sqrt{2\log(1/\delta)/N}$ we obtain, with probability
at least $1-\delta$,
\[
\|d_N\|_{L^2}^2
\le \frac{B_*^2}{N}\bigl(1+\sqrt{2\log(1/\delta)}\bigr)^2
\le \frac{B_*^2}{N}\bigl(2+4\log(1/\delta)\bigr).
\]
The right-hand side of \eqref{eq:fixed-population-bound} is independent
of $i$. Substituting this bound on $\|d_N\|_{L^2}^2$ therefore
gives the stated bound on the maximum.
\end{proof}

\begin{remark}
\label{rem:common-noise}
Suppose a Brownian motion $W^0$ is independent of the types and individual noises,
is observed by all players, and $c_{\theta_j}\,dW_t^0$ is added
to \eqref{eq:finite-game-state}, for a bounded continuous $c$.
Its second-order contribution
$\tfrac12\sum_{j,k}c_{\theta_j}c_{\theta_k}S_{jk}$
is independent of the state, leaving the Riccati equations
\eqref{eq:asy-Riccati} and \eqref{eq:finite-Riccati} unchanged.
The mean must then be conditioned on the common noise, as in
\cite{SiShi26}. The affine coefficient
$\eta_\theta$ satisfies a backward stochastic equation rather
than a terminal ODE evaluated along a realized mean path.
Extending Theorem~\ref{thm:main} requires an adapted consistency
system and conditional moment estimates. We leave this extension,
and state- or action-dependent common noise, for future work.
\end{remark}
% temperature.tex
\section{Temperature Dependence}
\label{sec:temperature}

The Gaussian policy mean is independent of $\gamma$, while
temperature enters the gain through the state covariance and the
variance-ratio term in \eqref{eq:exact-gain}. We vary the common
$\gamma$ with the type vector and mean input fixed.

\begin{proposition}
\label{prop:temperature}
For fixed $N,i,\theta^N$ and $m$,
there are nonnegative numbers
$a_{i,N}(\theta^N;m)$ and $b_{i,N}(\theta^N)$ such that
\begin{equation}
\label{eq:affine-temperature}
\varepsilon_{i,N}(\theta^N;m;\gamma)
=a_{i,N}(\theta^N;m)+\gamma b_{i,N}(\theta^N),\qquad\gamma>0.
\end{equation}
The slope is independent of $m$ and satisfies $b_{i,N}\le C/N$,
uniformly in the type vector. If $\bar Q=\bar H=0$ and
$\inf_{\theta\in\Theta}|D_\theta \barC_\theta|>0$, then
\begin{equation}
\label{eq:sharp-temperature}
c/N\le b_{i,N}(\theta^N)\le C/N\qquad(N\ge2),
\end{equation}
where $c,C>0$ depend only on the horizon, discount and coefficient bounds, including the stated
positive lower bounds.
\end{proposition}

At a mean field fixed point $m^*$ in the subclass of
Proposition~\ref{prop:temperature}, the lower slope bound and
$a_{i,N}\ge0$ give, after averaging over opponent types,
$c\gamma/N\le\mathcal E_{i,N}(\zeta;m^*)$.
Together with Theorem~\ref{thm:main}, this shows that the
$N^{-1}$ order is attained at each fixed positive temperature.

\begin{proof}
By \eqref{eq:finite-mean-ODE}, the mean $\nu$ and hence $f_j$ are
independent of $\gamma$. We separate the forcing proportional to
$\gamma$ in \eqref{eq:covariance-ODE}. The matrix $\Sigma^{(0)}$ solves
that equation at $\gamma=0$, and $\Sigma^{(1)}$ solves
\begin{equation}
\label{eq:covariance-parts}
\begin{aligned}
\dot\Sigma^{(1)}={}&\mathcal A^{\rm cl}\Sigma^{(1)}
+\Sigma^{(1)}(\mathcal A^{\rm cl})^\top\\
&+\operatorname{diag}_j\left\{(\ell_j^{\rm cl})^\top\Sigma^{(1)}\ell_j^{\rm cl}
+D_j^2/\widehat R_j\right\},\qquad\Sigma^{(1)}(0)=0.
\end{aligned}
\end{equation}
Uniqueness of the linear equation gives
\begin{equation}
\label{eq:covariance-temperature}
\Sigma_\gamma=\Sigma^{(0)}[m]+\gamma \Sigma^{(1)}.
\end{equation}
Since $\Sigma_\gamma$ is a covariance matrix for every $\gamma>0$,
the closedness of the positive semidefinite cone gives
\[
\Sigma^{(0)}=\lim_{\gamma\downarrow0}\Sigma_\gamma\ge0,
\qquad \Sigma^{(1)}=\lim_{\gamma\to\infty}\gamma^{-1}\Sigma_\gamma\ge0.
\]
The equation for $\Sigma^{(1)}$ contains neither $\nu$ nor $f_j$;
the feedback coefficients in it are determined by $P$ and are
independent of $m$. The same is true of $\widehat R_i^N$ and $\omega_i$.

Substitution of \eqref{eq:covariance-temperature} into
\eqref{eq:feedback-second-moment} and \eqref{eq:exact-gain}
gives the affine form, with
$a_{i,N}=\lim_{\gamma\downarrow0}\varepsilon_{i,N}\ge0$ and slope
\begin{equation}
\label{eq:temperature-slope}
b_{i,N}=\frac12\int_0^T e^{-\beta t}
\left[\widehat R_i^N\omega_i^\top \Sigma^{(1)}\omega_i
+\varrho_i-1-\log\varrho_i\right]dt.
\end{equation}
Positive semidefiniteness and
$\varrho_i-1-\log\varrho_i\ge0$ give $b_{i,N}\ge0$.

For the slope bound, apply the integral estimate leading to
\eqref{eq:covariance-bound} to the equation for $\Sigma^{(1)}$:
\[
\|\Sigma^{(1)}(t)\|_{\rm op}
\le C\int_0^t\|\Sigma^{(1)}(s)\|_{\rm op}ds+Ct.
\]
Gronwall gives $\sup_{t\le T}\|\Sigma^{(1)}(t)\|_{\rm op}\le C$,
independently of $N,m,x_0$. Equations
\eqref{eq:feedback-coefficient-bounds},
\eqref{eq:variance-ratio-bound} and \eqref{eq:temperature-slope}
therefore give
\[
0\le b_{i,N}\le C\int_0^T(|\omega_i|^2+N^{-2})dt\le C/N.
\]

For the lower bound, compactness of $\Theta$ and continuity of the
coefficients give
\[
\frac{D_j^2}{\widehat R_j(t)}
\ge\frac{\inf_\theta|D_\theta\barC_\theta|^2}
{\sup_\theta|\barC_\theta|^2\,\sup_{\theta,s}\widehat R_\theta(s)}
=:c_0>0.
\]
If
$\partial_t\Psi(t,s)=\mathcal A^{\rm cl}(t)\Psi(t,s)$ and
$\Psi(s,s)=I_N$, variation of constants gives
\begin{equation}
\label{eq:covariance-duhamel}
\Sigma^{(1)}(t)=\int_0^t\Psi(t,s)
\operatorname{diag}_j\left\{
(\ell_j^{\rm cl})^\top\Sigma^{(1)}(s)\ell_j^{\rm cl}
+D_j^2/\widehat R_j(s)\right\}\Psi(t,s)^\top ds.
\end{equation}
As $\|\Psi(t,s)^{-1}\|_{\rm op}\le e^{C(t-s)}$ and
$\Sigma^{(1)}\ge0$, this implies
\begin{equation}
\label{eq:covariance-lower}
v^\top\Sigma^{(1)}(t)v
\ge c_0\int_0^te^{-2C(t-s)}ds\,|v|^2
\ge c_T|v|^2,\qquad t\in[T/2,T].
\end{equation}
Here $c_T=c_0(T/2)e^{-2CT}>0$; the forcing produces positivity
despite zero initial covariance, cf. \cite[Lemma~3.7]{GRZ24}.

To use \eqref{eq:covariance-lower} in \eqref{eq:temperature-slope},
we next bound $|\omega_i|$ from below near $T$.
When $\bar Q=\bar H=0$, $S(T)=H_ie_ie_i^\top$.
For $j\ne i$,
\[
\dot S_{ij}=(\beta-A_i-A_j^{\rm cl})S_{ij}
-\frac{(S\barA^N)_i+(S\barA^N)_j}{N}
-\sum_kS_{kk}(\ell_k)_i(\ell_k)_j+\Lambda_i\Lambda_j/\widehat R_i^N.
\]
The row bounds give
$|(S\barA^N)_i|+N|(S\barA^N)_j|\le C$,
and $|\Lambda_i\Lambda_j|\le C/N$. In the diffusion sum,
the $k=i$ term has magnitude at most $CS_{ii}/N$, since
$(\ell_i)_i=C_i+\barC_i/N$ and $(\ell_i)_j=\barC_i/N$.
The $k=j$ term has magnitude at most $CS_{jj}/N$ and the remaining
terms are bounded in sum by $C\tr S/N^2$. Together with the diagonal
equation these bounds imply
$|\dot S_{ii}|\le C$ and $|\dot S_{ij}|\le C/N$.
Let $H_0=\inf_\theta H_\theta>0$ and
$d_0=\inf_\theta|D_\theta \barC_\theta|>0$.
Choose $\delta\in(0,T/2]$ sufficiently small that
$C\delta\le H_0/2$ and
$\sup_\theta|B_\theta|\,C\delta\le d_0H_0/4$ for the preceding
uniform derivative constants. Integrating backward from $T$ gives
\[
S_{ii}(t)\ge H_0/2,\qquad |S_{ij}(t)|\le C\delta/N,
\qquad t\in[T-\delta,T],\quad j\ne i.
\]
The triangle inequality then gives
\[
|B_iS_{ij}+D_iS_{ii}\barC_i/N|
\ge\frac{d_0H_0}{2N}-\frac{\sup_\theta|B_\theta|C\delta}{N}
\ge\frac{d_0H_0}{4N}.
\]
Since $\widehat R_i^N$ is uniformly bounded above,
$|(\omega_i)_j|\ge c/N$ for $j\ne i$ on this interval.
Therefore $|\omega_i|^2\ge c(N-1)/N^2\ge c/(2N)$ for $N\ge2$.
Since $[T-\delta,T]\subset[T/2,T]$, \eqref{eq:covariance-lower}
and the positive lower bound on $\widehat R_i^N$ give
\[
b_{i,N}\ge\frac12\int_{T-\delta}^T e^{-\beta t}\widehat R_i^N
\omega_i^\top\Sigma^{(1)}\omega_i\,dt\ge\frac{c\delta}{N}.
\]
The interval length $\delta>0$ is independent of $N$, which proves
the required lower bound.
\end{proof}

The leading exploration contribution in this subclass comes from
$\Sigma^{(1)}$ in the mean-feedback part of the gap. The variance-ratio
term is $O(\gamma/N^2)$.

In the strict game \eqref{eq:ordinary-state} without the entropy term,
completing the square in both Hamiltonians gives the same feedback means.
For strict deviations with a nonexplosive, pathwise unique strong solution
and finite expected action energy, the verification in
Proposition~\ref{prop:finite-response} applies with $a_t^2$ as the action
moment and $\widehat R_i^N(a_t-u_i^{*,N})^2/2$ as the residual.
The minimizing affine feedback is admissible by linear growth and
Lipschitz continuity of its closed-loop coefficients.
The strict candidate has the same consistency map. Its gain is $a_{i,N}$.
Using its covariance $\Sigma^{(0)}$ in the mean and covariance estimates
gives \eqref{eq:fixed-population-bound} at $\gamma=0$. Averaging over
opponent types gives
\[
\bE_i^\zeta[a_{i,N}(\theta^N;m)]
\le C_M\left(\frac1N+\|\mathscr T[m]-m\|_{L^2(0,T)}^2\right),
\qquad \|m\|_\infty\le M .
\]

When $D=0$, the covariance and gap formulas give
\[
\Sigma^{(1)}=0,\qquad\varrho_i=1,\qquad
\varepsilon_{i,N}(\theta^N;m;\gamma)=\varepsilon_{i,N}(\theta^N;m;1).
\]
Applying Theorem~\ref{thm:main} at $\gamma=1$ then gives a
$C/N$ Nash bound at a mean field fixed point, uniformly in $\gamma>0$.
% numerics.tex
\section{Numerical Illustrations}
\label{sec:numerics}

We evaluate \eqref{eq:exact-gain} at fixed types by solving the scalar
control, finite response and moment equations, using
\eqref{eq:feedback-second-moment} in place of Brownian sampling. Take
\begin{gather}
\Theta=[0,1],\; F=\operatorname{Unif}[0,1],\;
T=0.6,\;\beta=0.2,\;x_0=0.7,\; A_\theta=-0.4+0.2\theta,\notag\\
\barA_\theta=0.1+0.1\theta,\;B_\theta=0.8+0.4\theta,\;
D_\theta=0.65+0.15\theta,\;C_\theta=0.10+0.05\theta,\;
\barC_\theta=0.7+0.2\theta,\notag\\
\sigma_\theta=0.5+0.25\theta,\;R_\theta=1+0.3\theta,\;
Q_\theta=1+0.4\theta,\;H_\theta=1+0.2\theta.
\label{eq:numerical-parameters}
\end{gather}

Comparison in the Riccati equation gives $0.7\le P_\theta\le H_\theta\le1.2$.
The bounds $\|A^{\rm cl}\|_\infty<1.13$, $\|b^\mu\|_\infty<0.31$,
$\|b^\eta\|_\infty<0.23$ and $\|B^2/\widehat R\|_\infty<0.84$
then give $L_T<0.75$ in \eqref{eq:contraction-bound}.
Picard iteration from $m^{(0)}=0$ with 64-point Gauss--Legendre
quadrature and 801 time nodes gives a final grid residual of
$1.1\times10^{-14}$.

\begin{figure}[!htbp]
\centering
\includegraphics[width=0.90\textwidth]{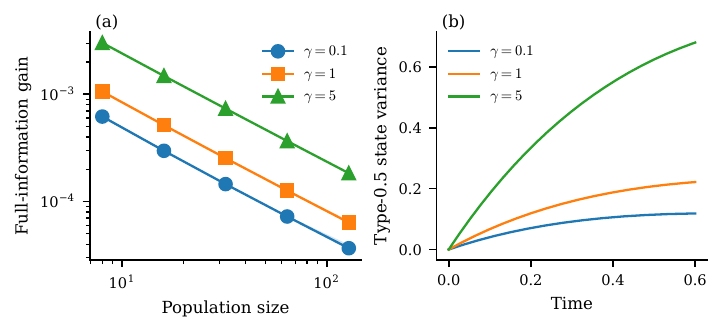}
\caption{(a) Mean full-information gain over six opponent-type populations,
with shaded sample ranges and own type $0.5$.
(b) Prescribed-mean state variance at type $0.5$.
The feedback mean is the same at all temperatures}
\label{fig:gain-state}
\end{figure}

\begin{figure}[H]
\centering
\includegraphics[width=0.90\textwidth]{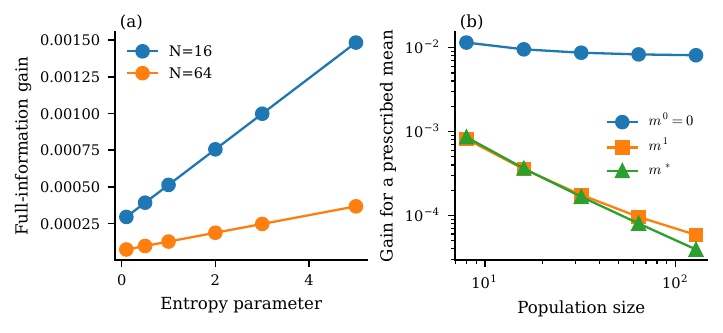}
\caption{(a) Gain versus temperature for fixed type populations, with
$\bar Q=\bar H=0$. (b) Gain at $\gamma=1$ for three mean inputs in midpoint
type populations, with $\bar Q_\theta=0.4+0.2\theta$ and
$\bar H_\theta=0.3+0.1\theta$}
\label{fig:temperature-residual}
\end{figure}

Except in Figure~\ref{fig:temperature-residual}(b), take
$\bar Q=\bar H=0$. For Figure~\ref{fig:gain-state}(a),
at each $N\in\{8,16,32,64,128\}$,
fix player~1's type at $0.5$ and draw six independent
opponent-type populations. Use the same populations at
$\gamma=0.1,1,5$. The plot shows their
mean gains and sample ranges. The fitted log--log slope at
$\gamma=1$, over $N=16,32,64,128$, is $-1.005$, consistent
with the $N^{-1}$ gain bounds at a mean field fixed point
(Theorem~\ref{thm:main} and Proposition~\ref{prop:temperature}).
Panel (b) shows the variance of \eqref{eq:limiting-state} at type $0.5$
under the computed mean.

In Figure~\ref{fig:temperature-residual}(a), the scaled slopes
$Nb_{1,N}$ are $3.879\times10^{-3}$ at $N=16$ and $3.841\times10^{-3}$
at $N=64$, consistent with \eqref{eq:sharp-temperature}. At $N=16$,
replacing the consistent input by $m=0$ changes the slope by less
than $5\times10^{-14}$. Panel (b) uses midpoint types
$\theta_j=(j-1/2)/N$. Its mean-state cost weights affect only
$\kappa_\theta$ in \eqref{eq:kappa}. One update
$m^{(1)}=\mathscr T[0]$ reduces the gain from using $m^{(0)}=0$,
and the gain decreases with $N$ at the computed consistent input.
This behavior agrees with \eqref{eq:fixed-population-bound}.

\begin{samepage}
We use DOP853 with tolerances $2\times10^{-10}$ (relative) and
$2\times10^{-12}$ (absolute), and 1601 gain integration nodes. At
$\gamma=1$, direct cost minus response value agrees with
\eqref{eq:exact-gain} within $2\times10^{-11}$. For the $N=16$ population
in Figure~\ref{fig:temperature-residual}(a) at $\gamma=1$, using
128 quadrature, 1601 mean and 3201 integration nodes and a tighter relative
tolerance changes the gain by less than $10^{-8}$ relative to its refined value.

\end{samepage}
\section{Conclusions}
\label{sec:conclusions}
At a mean field fixed point and fixed temperature, the gain averaged
over opponent types decays at least as fast as the reciprocal of the
population size. The affine temperature formula isolates the additional
gain due to exploration. This contribution attains the same rate when
the mean-state costs vanish and both diffusion loadings are bounded
away from zero.

% declarations.tex
\section*{Declarations}

\paragraph{Funding.} The author received no research funding for this work.

\paragraph{Conflict of Interest.}
The author has no relevant financial or non-financial interests to disclose.

\paragraph{Data and Code Availability.}
The code and numerical output are available from the corresponding
author on reasonable request.

\paragraph{Use of Generative AI.}
ChatGPT assisted with mathematical exploration, literature review,
numerical code, and drafting and revision. DeepSeek provided feedback
on the manuscript, and Rethlas assisted with mathematical checks.
The author takes responsibility for the final manuscript and results.

\providecommand{\doi}[1]{\href{https://doi.org/#1}{\nolinkurl{doi:#1}}}
% References are included here so the preprint can be built directly with pdfLaTeX.

\end{document}